\documentclass[11pt]{amsart}

\usepackage{bbm}
\usepackage{amssymb}
\usepackage{graphicx}
\usepackage{enumerate}
\usepackage{xcolor}
\usepackage{mathabx}
\usepackage {varioref} 
\usepackage{tikz}

\newtheorem{theorem}{Theorem}[section]
\newtheorem{proposition}[theorem]{Proposition}
\newtheorem{lemma}[theorem]{Lemma}

\RequirePackage{coloralias}
\definecoloralias {utahheadingcolor}{utahred}
\newcommand {\figlabel}    [1] {\label{fig:#1}}
\newcommand {\figref}      [1] {{\color {utahheadingcolor}\textbf{Figure~\ref{fig:#1}}} \vpageref{fig:#1}}

\theoremstyle{definition}

\newenvironment{psmallmatrix}
  {\left(\begin{smallmatrix}}
  {\end{smallmatrix}\right)}

\newcommand{\appsection}[1]{\let\oldthesection\thesection
\renewcommand{\thesection}{Appendix \oldthesection}
\section{#1}\let\thesection\oldthesection}

\theoremstyle{remark}

\newtheorem{remark}[theorem]{Remark}

\DeclareMathOperator{\SL}{\textbf{SL}}
\DeclareMathOperator{\GL}{\textbf{GL}}
\DeclareMathOperator{\Stab}{\textbf{Stab}}
\DeclareMathOperator{\Lk}{Lk}
\DeclareMathOperator{\St}{St}
\DeclareMathOperator{\Span}{Span}

\DeclareMathOperator{\im}{im}

\def\Z{{\mathbb{Z}}}
\def\Q{{\mathbb{Q}}}
\def\C{{\mathcal{C}}}
\def\A{{\mathcal{A}}}
\def\S{{\mathcal{S}}}

\def\G{{\textbf{G}}}

\newcommand{\Proof}{{\it Proof. }}

\newcommand{\Sect}{\mathfrak{S}}
\newcommand{\Vect}[2]{\begin{psmallmatrix} #1\\#2\end{psmallmatrix}}
\newcommand{\Congk}{\Lambda_k}
\newcommand{\CH}{Ch}
\newcommand{\Eta}{H}

\newcommand{\SecRef}[1]{Section~\ref{sec:#1}}

\newcommand{\SecsRef}[2]{Sections~\ref{sec:#1} and~\ref{sec:#2}}
\newcommand{\sectionChapter}{section }

\begin{document}
\bibliographystyle{amsplain}
\title{An Infinite Dimensional Virtual Cohomology Group of $\SL_3(\Z[t])$}
\author{\textrm{Matthew Goroff}}


\begin{abstract}
We prove that $\SL_3(\Z[t])$ has a finite index subgroup $\Gamma$ such that $H^2(\Gamma; \Q)$ is infinite dimensional. The proof uses the geometry of the Euclidean building for $\SL_3(\Q((t^{-1})))$.
\end{abstract}

\maketitle

\section{Introduction} \label{s0}

In \cite{Sus}, Suslin proves that $\SL_3(\Z[t])$ is finitely generated by elementary matrices. In \cite{KM}, Kristic and McCool prove that $\SL_3(\Z[t])$ is not finitely presented. In \cite{BMW}, Bux, Mohammadi, and Wortman prove that $\SL_3(\Z[t])$ is not $FP_2$.

In this paper, we use and add to these results by showing that $\SL_3(\Z[t])$ has a finite index subgroup with an infinite dimensional second cohomology group with coefficients in $\Q$. Specifically, consider the homomorphism $\iota\colon\Z[t]\to\Z/3\Z$ where $t\mapsto0$ and integers are reduced mod 3. This induces a homomorphism $\hat{\iota}\colon\SL_3(\Z[t])\to\SL_3(\Z/3\Z)$. Let $\Gamma = \ker(\hat{\iota})$. We will prove the following:

\noindent \textbf{Theorem 1. } {\it $H^2(\Gamma; \Q)$ is infinite dimensional.}

The structure of this paper borrows from \cite{CK}, where Cesa and Kelly attempt to prove that $H^2(\SL_3(\Z[t]);\Q)$ is infinite dimensional. In \cite{CK}, they define an infinite family of cocycles, but they do not prove that these cocycles are non-zero. For an explanation of why this is not done in \cite{CK}, see remark~\ref{rem:rebuttal}. We manage to avoid this by restricting to the finite index subgroup $\Gamma$. We also take a closer look at the spherical building for $\SL_3\Q$ as a way to construct cycles which are evaluated nontrivially. In \cite{CK}, they borrow a Morse function defined in \cite{BKW}, but it is used on a subspace of the building which does not necessarily have connected descending links. We avoid this problem by not using a Morse function, and instead find an alternate way to extend local disks to global disks in the building. This process is done in \SecRef{globalCocycles} and takes up a substantial portion of the paper.

\subsection{Outline of Paper}
The structure of this paper is modeled on \cite{Wort}, where Wortman proves a similar result for arithmetic groups over function fields. The main difference is that much effort is spent in this paper on extending local disks about vertices to global disks whose boundaries are contained in a neighbordhood of a $\Gamma$-orbit in a Euclidean building, whereas Wortman utilizes a Morse function developed in \cite{BKW}. We also face the challenge that the Euclidean building we are working with is not locally compact.

In \SecRef{preliminaries}, we review some key notions about the Euclidean building for $\SL_3(\Q((t^{-1})))$ and the spherical building for $\SL_3\Q$ that will be used in this paper. We briefly explore the connection between the two, as well as the spherical building which is the visual boundary of the Euclidean building. In addition, we establish notation that will be used to denote upper triangular subgroups throughout the paper. Then we illustrate, via a lemma, why we chose $\Gamma$ in such a way.

In \SecRef{localStars}, we define several key maps which allow us to translate the action of a stabilizer of a vertex on that vertex's star to the action of an upper triangular subgroup on a spherical building. Then in \SecRef{sequenceOfVertices}, we choose a specific sequence of vertices which will be used throughout the paper.

In \SecRef{connectedComplex}, we construct a 2-connected space for $\Gamma$ to act on, and show that it has a large contractible subspace.

In \SecRef{localCocycles}, we shift our attention to the spherical building for $\SL_3\Q$ and define a nontrivial cocycle on a quotient space of it. Later in the paper, we extend this cocycle, and the cycle which it evaluates to be non-zero, to a cocycle and cycle on the larger building.

In \SecsRef{congruenceSubgroups}{reducingDisks}, we define a series of unipotent subgroups in $\SL_3(\Q[t])$ and show that after modding out by them, we can bound the height of certain disks in the building.

In \SecRef{globalCocycles}, we define an infinite sequence of cocycles modeled on the local cocycles found in \SecRef{localCocycles}. Then we extend the local cycles found in \SecRef{localCocycles} to global disks in the Euclidean building and cycles in the quotient space.

Finally, in \SecRef{mainProof}, we prove Theorem 1 using a spectral sequence argument. 

\section{Preliminaries}
\label{sec:preliminaries}

We will make heavy use of Euclidean and spherical buildings in this paper. We will state some key notions below, but for reference, see \cite{Build}.

Let $X$ be the Euclidean building for the group $\SL_3(\Q((t^{-1})))$. We let $\A$ be the apartment in $X$ corresponding to the diagonal subgroup of $\SL_3(\Q((t^{-1})))$ and $x_0$ be the vertex fixed by $\SL_3(\Q[[t^{-1}]])$. Recall that there is an action of $\GL_3(\Q((t^{-1})))$ on $X$ which is transitive on vertices. For $s_1,s_2,s_3\in\Q((t^{-1}))^\times$, we let $D(s_1,s_2,s_3)$ be the diagonal matrix with the $s_i$ along the diagonal. Let $\C$ be the chamber with vertices $x_0$, $D(t,1,1)x_0$, and $D(t,t,1)x_0$. Let $\Sect$ be the sector in $\A$ with vertex $x_0$ and containing $\C$. We recall that $\Sect$ is a fundamental domain for the action of $\SL_3(\Q[t])$ on $X$ \cite{Sou}. The vertices of the sector $\Sect$ can be given uniquely by diagonal matrices of the form $D(t^{m_1},t^{m_{2}},1)$ where $m_1\ge m_2\ge0$. That is, each vertex of $\Sect$ is a translate of $x_0$ by a unique such diagonal matrix.

We will call $x_0$ the standard vertex in $X$, $\C$ the standard chamber in $X$, and $\A$ the standard apartment in $X$.

Notice that $\Gamma < \SL_3(\Z[t]) < \SL_3(\Q[t]) < \SL_3(\Q((t^{-1})))$, so each has a well-defined action on $X$. For a cell $\sigma < X$, we write $K_\sigma = \Stab_{K}(\sigma)$ where $K$ is any group which acts on $X$.

Let $\Xi$ be the spherical building for $\SL_3(\Q)$. $\Xi$ can be realized as the flag complex of subspaces of $\Q^3$ where a basis for $\Q^3$ will give an apartment in $\Xi$. Let $e_i$ be the $i$th standard coordinate vector in $\Q^3$, then the basis $\{e_1,e_2,e_3\}$ gives an apartment in $\Xi$ which we will call $\A_\Xi$. Equivalently, we can describe simplices in $\Xi$ by their stabilizers in $\SL_3(\Q)$, which will be parabolic subgroups. We will call $\A_\Xi$ the standard apartment in $\Xi$, and the chamber which corresponds to the complete flag $\langle e_1\rangle < \langle e_1,e_2\rangle < \Q^3$ the standard chamber.

As $\Q$ is the residue field for $\Q((t^{-1}))$, the link about any vertex in $X$ is naturally identified with $\Xi$. In \SecRef{localStars}, we will detail this relationship more closely.

We let $X(\infty)$ be the visual boundary of $X$. This can be naturally identified with the spherical building for $\SL_3(\Q((t^{-1})))$, so we can describe simplices in $X(\infty)$ by partial flags in $(\Q((t^{-1})))^3$. We will write $\A(\infty)$ to denote the visual boundary of the apartment $\A < X$. This is the apartment of $X(\infty)$ given by the basis $\{e_1,e_2,e_3\}$ of $(\Q((t^{-1})))^3$.

We will reference upper triangular subgroups frequently in this paper, so adopt the following notation.
\[
U(K) \text{ is the subgroup of upper triangular matrices in $K$ with 1's on the diagonal}
\]
\[
B(K) \text{ is the subgroup of upper triangular matrices in $K$}
\]
where $K$ can be any subgroup of $\SL_3(\Q((t^{-1})))$ which we reference in the paper. It is always the case that $U(K) \subseteq B(K)$, but when $K=\Gamma$, they are equal, which we prove below and is the motivating reason we chose $\Gamma$. 
\begin{lemma}
$B(\Gamma) = U(\Gamma)$.
\end{lemma}
\Proof Notice that if $u\in B(\Gamma)$, then the product along the diagonal must be 1, so each element must be either 1 or -1 (the only units of $\Z)$. But each must be taken to 1 under $\iota$, that is, each must be 1 and $u\in U(\Gamma)$.\qed

\section{Stars about Vertices in $\Sect$}
\label{sec:localstars}
\label{sec:localStars}

In this \sectionChapter, we give an explicit map from $\St(y)$ to $\Xi$ for any vertex $y$ in the interior of $\Sect$, which we denote $\mathring{\Sect}$. In addition, we provide an equivariant map from $\Stab_{\SL_3(\Q[t])}(y)$ to $\SL_3\Q$. Throughout this paper, we use open cells and open stars.

First we note that there is a well-defined map from the chambers in $\St(y)$ to $\Lk(y)$:
\begin{eqnarray*}
\downarrow_y\colon\CH(\St(y))&\to& \Lk(y)\\
\tau &\mapsto& \widebar{\tau}\cap\Lk(y)
\end{eqnarray*}
So in order to define a map from the chambers of $\St(y)$ to the spherical building, it is sufficient to define a map from $\Lk(y)$.

Recall that each vertex in the standard sector can be given uniquely by a certain diagonal matrix in $\GL_3(\Q((t^{-1})))$. For each $y\in\Sect$, let $D_y$ be this diagonal matrix, that is, $y = D_y\cdot x_0$.

A simple valuation argument gives us that $\Stab_{\SL_3(\Q[t])}(y)\subseteq\hat{U}(\SL_3(\Q[t]))$. Let $u\in\Stab_{\SL_3(\Q[t])}(y)$, then $D_y^{-1} u D_y\in\Stab(x_0) = \SL_3(\Q[[t^{-1}]])$. Consider the homomorphism from $\Q[[t^{-1}]]\to\Q$ which is the identity on $\Q$ and takes $t^{-1}$ to 0; this induces a homomorphism from $\SL_3(\Q[[t^{-1}]])\to\SL_3\Q$, let's call this homomorphism restriction to $\SL_3\Q$ (so we write $|_{\SL_3\Q}$). This restriction homomorphism takes upper triangular matrices to upper triangular matrices. So we can define
\begin{eqnarray*}
\rho_y\colon\Stab_{\SL_3(\Q[t])}(y)&\to&B(\SL_3\Q)\\
u&\mapsto& (D_y^{-1} u D_y)|_{\SL_3\Q}
\end{eqnarray*}

As we stated in the prior \sectionChapter, the link of $x_0$ is exactly $\Xi$. If we identify $\Lk(x_0)\cap\C$ with the standard chamber in $\Xi$, then we can push around this identification with $\SL_3\Q$. As these two spaces are identified and have an identical action by $\SL_3\Q$, we will refer to them interchangeably. 

Now we define a map from the chambers of $\St(y)$ to the spherical building:
\begin{eqnarray*}
\eta_y\colon\CH(\St(y)) &\to& \Xi\\
\tau &\mapsto& D_y^{-1}(\downarrow_y(\tau))
\end{eqnarray*}

\begin{lemma}
For each $y\in\mathring{\Sect}$, $\rho_y$ and $\eta_y$ are equivariant. That is, for $E\in\CH(\St(y))$ and $u\in\Stab_{\SL_3(\Q[t])}(y)$:
\[
\eta_y(u\cdot E) = \rho_y(u)\cdot \eta_y(E)
\]
\end{lemma}
\Proof First note that for any $q\in\SL_3(\Q[[t^{-1}]])$ and $x\in\Lk(x_0)$, $q\cdot x = q|_{\SL_3\Q}\cdot x$. So we can manipulate:
\[
\eta_y(u\cdot E) = D_y^{-1}(\downarrow_y(u\cdot E))) = D_y^{-1} u D_y D_y^{-1} (\downarrow_y (E)) = (D_y^{-1} u D_y)\cdot \eta_y(E) = \rho_y(u)\cdot \eta_y(E)
\]\qed

In addition, we prove the following two lemmas, which will be useful later:

\begin{lemma}
\label{lem:commutators}
For any $y\in\mathring{\Sect}$, $\rho_y(U(\Stab_{\SL_3(\Q[t])}(y)))$ surjects onto $U(\SL_3\Q)$. In addition, $\rho_y([U(\Stab_{\SL_3(\Q[t])}(y)), U(\Stab_{\SL_3(\Q[t])}(y))])$ surjects onto $[U(\SL_3\Q), U(\SL_3\Q)]$.
\end{lemma}
\Proof The first claim is a consequence of the restriction map $|_{\SL_3\Q}$ being the identity on $\SL_3\Q$. So for any $h\in U(\SL_3\Q)$, we have that $D_yhD_y^{-1}\in\Stab_{\SL_3(\Q[t])}(y)$ and $\rho(D_yhD_y^{-1}) = h$. The same logic gives us the second claim.\qed

\begin{lemma}
\label{lem:integerImage}
For any $y\in\mathring{\Sect}$, $\rho_y(U(\Stab_{\SL_3(\Z[t])}(y)))$ surjects onto $U(\SL_3\Z)$.
\end{lemma}
\Proof Let $u\in U(\SL_3\Z)$. Then $D_yuD_y^{-1}\in U(\Stab_{\SL_3(\Z[t])}(y))$, and $\rho_y(D_yuD_y^{-1}) = u$.\qed

\section{A Sequence of Vertices in $\Sect$}
\label{sec:sequenceOfVertices}

In this \sectionChapter{}, we describe a sequence of vertices in the sector $\Sect$ which will be used throughout the paper.

In \cite{Sch}, Schulz gives a description of open hemisphere complexes in spherical buildings which we will describe here and use for this paper. Any spherical building is endowed with the angular metric $d_\Xi$. For any vertex $x\in\Xi$, the open hemisphere complex opposite $x$ is the full subcomplex spanned by $\{y\in\Xi^{(0)}\ |\ d_\Xi(x,y) > \pi/2\}$. In \cite{Sch}, it is proved that such open hemisphere complexes are spherical of dimension $1$. Although the name is confusing, these open hemisphere complexes are in fact closed simplicial complexes.

Let $P\in\Xi$ be the parabolic subgroup of $\SL_3\Q$ which stabilizes the subspace $\langle e_1,e_2\rangle < \Q^3$. We will use $P$ to denote both the subgroup and the corresponding vertex in $\Xi$. Let $\Omega<\Xi$ be the open hemisphere complex opposite $P$. For any vertex $y\in\Sect$, we let $\Omega(y)=\eta_y^{-1}(\Omega)<\St(y)$. 

\begin{lemma}
$\Omega = P\cdot \widebar{\St}(\langle e_3\rangle)$, where we are taking the closed star about the vertex $\langle e_3\rangle\in\Xi$.
\end{lemma}
\Proof First we note that $P\cdot\A_\Xi = \Xi$ because $\Xi$ is a building. If we focus our attention on the apartment $\A_\Xi$, we can see that each vertex is connected by an edge to either $\langle e_1,e_2\rangle = P$ or $\langle e_3\rangle$. As $\Xi$ is irreducible, each chamber has diameter less than $\pi/2$. So the vertices in $\A_\Xi$ which are in $\Omega$ are exactly those that are connected to $\langle e_3\rangle$ by an edge. The span of these is precisely the closed star of $\langle e_3\rangle$ in $\A_\Xi$. As $P$ acts by isometries on $\Xi$, we can push this closed star around by $P$ to obtain all of $\Omega$ while also showing that every other vertex is within $\pi/2$ of $P$.\qed\\

For any $y\in\Sect$, the chambers of $\Omega(y)\cap\A$ can be extended to sectors in the apartment $\A$ of $X$. For any such chamber $C$, let $\widehat{C}$ denote this sector intersected with $\Sect$.

\begin{lemma}
For any $C<\Omega(y)\cap\A$, $\widehat{C}$ is compact.
\end{lemma}
\Proof The sector of $\A$ associated to $C$ is the union of every ray in $\A$ starting at $y$ and going through $C$. Chambers in irreducible buildings have radius at most $\pi/2$, so the half-space of $\A$ composed of every ray starting at $y$ with an angle greater than $\pi/2$ from $P$ intersects $\Sect$ in a compact set. Each $\widehat{C}$ is contained in here, and hence is also compact.\qed\\

Now let
\[
\widehat{\Omega}_\A(y) = \bigcup_{C < \Omega(y)\cap\A} \widehat{C}
\]
which by construction is a compact subset of $\Sect$.

\begin{lemma}
There exists $y_1,y_2,\dots\in\Sect$ such that for all $k\ge1$:
\begin{enumerate}[i)]
\item $y_k$ is not in the closed star of $\partial\Sect$, and
\item for all $i<k$, $y_k\notin\widehat{\Omega}_\A(P)(y_i)$
\end{enumerate}
\end{lemma}
\Proof Let $y_1\in\Sect$ be a vertex which is not in the closed star of the boundary of $\Sect$. Then we define $y_k$ recursively for $k>1$ to be a vertex which is neither in the closed star of the boundary of $\Sect$ nor in $\widehat{\Omega}_\A(P)(y_i)$ for all $i<k$. As each $\widehat{\Omega}_\A(P)(y_i)$ is compact, and $\Sect$ minus the closed star of the boundary is an unbounded space, we can always pick another $y_k$.\qed

\section{Construction of a $2$-connected $\Gamma$-complex}
\label{sec:connectedComplex}

In this \sectionChapter{}, we construct a 2-connected $\Gamma$-complex $Y$, as well as a $\Gamma$-equivariant map $\psi\colon Y\to X$.

Let $X_0=\SL_3(\Z[t])\cdot\C$, and $\psi\colon X_0\to X$ be the inclusion map. We will add to $X_0$ in order to construct $Y$, and extend $\psi$ at the same time.

A result of Suslin \cite{Sus} tells us that $\SL_3(\Z[t])$ is generated by elementary matrices which each fix a point in the closed star about $x_0$ in $\A$, a subspace of $X_0$, hence $X_0$ is connected.

Now, for any cellular map $\sigma\colon S_{1}\to X_{0}$ and $\gamma\in\Gamma$, glue a $2$-disk to $X_0$ along the attaching map $\gamma\cdot\sigma$. And have $\psi$ take this disk to the unique filling disk of $\psi(\gamma\cdot \sigma(S_{1}))$ in $X$. Call the resulting space $X_{1}$.

Now repeat the process an additional time, glueing in $3$-disks and call the resulting space $Y$. Let $\psi$ take each 3-disk to $X$ in a continuous fashion that is $\Gamma$-invariant.

By construction, this space has to be 2-connected, as we have glued in any required filling disks and filling spheres. Importantly, notice that the $1$-skeleton of $Y$ is the 1-skeleton of $X_0=\SL_3(\Z[t])\cdot\C$.

\subsection{$\Gamma\setminus X_0$ is Contractible}
Ultimately, we will be working with the quotient space $\Gamma\setminus Y$, but in this section, we prove the following convenient fact:
\begin{proposition}
\label{prop:contractibility}
The quotient space $\Gamma\setminus X_0$ is contractible.
\end{proposition}
We do this with the several following lemmas. First we remark that the action of $\SL_3(\Q((t^{-1})))$ is type preserving on simplices of $X$. So when we write that two simplices are of the same type in $X_0$, we mean that there exists some $h\in\SL_3(\Z[t])$ which takes one to the other. For any simplex $\S\subseteq  X_0$ and vertex $x\in\S$, let
\[
\S\setminus x = \Span\{v\in\S^{(0)}\ |\ v\neq x\}
\]
that is, the span of the rest of the vertices of $\S$. If $\S=x$, then this is the empty set.
\begin{lemma}
\label{lem:inStandard}
Let $\S\subseteq\C$ be a simplex which includes $x_0$. Let $\S'\subseteq X_0$ be a simplex of the same type such that $\S'\cap\S = \S\setminus x_0$. Then there exists $\gamma\in\Gamma$ such that $\gamma\S'=\S$.
\end{lemma}
\Proof First we address the case where $\S = x_0$. Thus $\S'=hx_0$ for some $h\in\SL_3(\Z[t])$. Let $h|_\Z$ be the element in $\SL_3\Z$ where we restrict each entry to the constant term. This gives an element of $\SL_3\Z$, and furthermore, $\hat{\iota}(h)=\hat{\iota}(h|_\Z)$. So consider the element $h(h|_\Z)^{-1}\in\SL_3(\Z[t])$. We have that $h(h|_\Z)^{-1}\cdot x_0 = h\cdot x_0$ as $\SL_3\Z$ fixes $x_0$. Also, $h(h|_\Z)^{-1}$ is in the kernel of $\hat{\iota}$ and thus is an element of $\Gamma$. So we have provided an element of $\Gamma$ which takes $\S$ to $\S'$.


Now we address the case where $\S\setminus x_0$ is a vertex. There are only two such choices for $\S$: when $\S\setminus x_0= D(t,1,1)x_0$ and $\S\setminus x_0=D(t,t,1)x_0$. Either way, there exists $h\in\SL_3(\Z[t])$ such that $h\S = \S'$, and $h\in\Stab_{\SL_3(\Z[t])}(\S\setminus x_0)$.  If we examine $h|_\Z$, we notice that in the first case,
\[
h|_\Z = \begin{pmatrix}
\pm 1 & * & *\\
0 & * & *\\
0 & * & *
\end{pmatrix}\in\SL_3(\Z)
\]
which also fixes $D(t,1,1)x_0 = \S\setminus x_0$. In the second case,
\[
h|_\Z = \begin{pmatrix}
* & * & *\\
* & * & *\\
0 & 0 & \pm1
\end{pmatrix}\in\SL_3(\Z)
\]
which fixes $D(t,t,1)x_0 = \S\setminus x_0$. So we can run a similar argument as in the above paragraph and see that $h(h|_\Z)^{-1}(\S) = \S'$ and $h(h|_\Z)^{-1}$ is in $\Gamma$ as it is in the kernel of $\hat{\iota}$.

If $\S\setminus x_0$ is an edge, then it is the edge that spans the two vertices in the paragraph above. We can run the same argument, and see that
\[
h|_\Z = \begin{pmatrix}
* & * & *\\
0 & * & *\\
0 & 0 &*
\end{pmatrix}\in\SL_3(\Z)
\]
which fixes the edge and $x_0$. So $h(h|_\Z)^{-1}\in\Gamma$ is our desired element.\qed\\

\begin{lemma}
Let $\S\subseteq X_0$ be a simplex which contains $x_0$ and $\S'\subseteq X_0$ be a simplex of the same type such that $\S\cap\S'=\S\setminus x_0$. Then $\Gamma\S = \Gamma\S'$.
\end{lemma}
\Proof There must exist some $r\in\SL_3(\Z[t])$ such that $r\S\subseteq \C$. Then $r\S'\cap r\S = r\S\setminus x_0$. So by lemma~\ref{lem:inStandard}, there exists a $\gamma\in\Gamma$ such that $\gamma r\S'=r\S$. Thus $r^{-1}\gamma r\S' = \S$. $\Gamma$ is normal in $\SL_3(\Z[t])$ and hence we have shown that $\Gamma\S=\Gamma\S'$.\qed
\begin{lemma}
\label{lem:boundariesMatch}
Let $\S_1$ and $\S_2$ be simplices in $X_0$ such that $\Gamma(\partial\S_1) = \Gamma(\partial\S_2)$ and $\Gamma x_0\in \Gamma(\partial\S_1)$. Then $\Gamma\S_1=\Gamma\S_2$.
\end{lemma}
\Proof If $\S_1$ is a translate of $x_0$, then we are done, so assume that it is not. By assumption and lemma~\ref{lem:inStandard}, there exists $\omega_i\in\Gamma$ such that $\omega_i x_0\in\S_i$. $\S_1\setminus(\omega_1 x_0)$ and $\S_2\setminus(\omega_2 x_0)$ are boundary components of the same type, and so by assumption, there exists some $\gamma\in\Gamma$ such that $\gamma( \S_1\setminus(\omega_1 x_0)) = \S_2\setminus(\omega_2 x_0)$. So $\omega^{-1}_2\gamma\S_1 \cap\omega_2^{-1}\S_2 = \omega_2^{-1}\S_2\setminus x_0$. Then lemma 5.3 tells us that there exists an $\delta\in\Gamma$ such that $\delta\omega_2\gamma\S_1 = \omega_2\S_2$. Rearranging, we find that:
\[
\omega_2^{-1}\delta\omega_2\gamma\S_1=\S_2
\] 
That is, that $\Gamma\S_1 = \Gamma\S_2$.\qed\\

{\it Proof of Proposition~\ref{prop:contractibility}.} 
Let $\Gamma p \in \Gamma\setminus X_0$ and let $\Gamma\S$ be a minimal simplex which carries $\Gamma p$ and $\Gamma x_0$. As $\Gamma x_0$ is a vertex in each chamber of $\Gamma \setminus X_0$, such a simplex always exists. Lemma~\ref{lem:boundariesMatch} above confirms that this minimal simplex is in fact unique. So take the deformation retraction on $\Gamma\S$ which pushes points down towards $\Gamma x_0$. To see that this is well-defined, we note that this description will always agree on the intersection of two simplices, and that $\Gamma\S$ is well-defined for any point $p$.\qed\\

\section{Local Cocycles}
\label{sec:localCocycles}

In this \sectionChapter{}, we restrict our view completely to $\Xi$ and the action of $U(\SL_3\Q)$ on $\Xi$. With that in mind, in this \sectionChapter{}, we let $U=U(\SL_3\Q)$.

Let $\mathfrak{C}' < \Xi$ be the chamber opposite the standard chamber in the standard apartment. That is, $\mathfrak{C}'$ corresponds to the parabolic subgroup which stabilizes the flag $\langle e_3\rangle < \langle e_2,e_3\rangle < \Q^3$. Now let $\Theta = [U,U]\setminus(U\cdot\mathfrak{C}')$. That is, take all translates of $\mathfrak{C}'$ and mod out by the commutator subgroup of $U$. Let $\pi\colon U\cdot\mathfrak{C}'\to\Theta$ be projection. We let $\mathfrak{C}$ be the equivalence class $\pi(\mathfrak{C}')=[U,U]\cdot\mathfrak{C}'$, a 1-cell in $\Theta$. There is an action of $U$ on $\Theta$: let $u ([U,U]u'\mathfrak{C}') = [U,U] uu'\mathfrak{C}'$. This action is equivariant with respect to the action of $U$ on $U\cdot\mathfrak{C}'$. 

The stabilizer of $\mathfrak{C}'$ in $\SL_3\Q$ is all lower triangular matrices, making $U_{\mathfrak{C}'}$ trivial. Thus, the stabilizer of $\mathfrak{C}$ in $U$ is exactly $[U,U]$. That is, the group $[U,U]\setminus U$ is in correspondence with the 1-cells of $\Theta$.

\begin{lemma}
$[U,U]\setminus U \cong \Q^2$.
\end{lemma}
\Proof Consider the map
\begin{eqnarray*}
\varepsilon\colon U &\to& \Q^2\\
u &\mapsto& (u_{1,2}, u_{2,3})
\end{eqnarray*}
where $u_{i,j}$ is the entry of $u$ in the $i$th row and $j$th column. A quick check verifies that $\varepsilon$ is a surjective homomorphism. The kernel of $\varepsilon$ is $[U,U]$, giving us an isomorphism between $[U,U]\setminus U$ and $\Q^2$.\qed

As we have a bijection between $\Q^2$ and the 1-cells of $\Theta$, we will use them interchangeably.

\begin{lemma}
Let $g,h\in \Q^2$ such that they differ in exactly one position. Then $g\cdot \mathfrak{C}$ and $h\cdot \mathfrak{C}$ are distinct 1-cells in $\Theta$ which share a vertex. 
\end{lemma}
\Proof First we note that the elementary matrices (or folds)
\[
e_{i,i+1}(k)\in U
\]
each fix exactly one vertex of $\mathfrak{C}$ when $k\neq 0$. So it suffices to show that there exists an elementary fold such that $(g\cdot\varepsilon(e_{i,i+1}(k)))\cdot \mathfrak{C}=h\cdot \mathfrak{C}$.

Let $\widebar{\varepsilon}\colon \Q^2\to U$ be the map
\[
g\mapsto u
\]
where $u_{i,i+1}=g_i$ and $u_{i,j}=0$ for $i<j+1$. Notice that $\varepsilon\circ\widebar{\varepsilon}=id$.\\
Assume that $g,h\in \Q^2$ differ in the $i$th position. Then the following precompositon by an elementary fold shows what we are looking for.
\[
g\cdot\varepsilon(e_{i,i+1}(h_i-g_i)))=\varepsilon(\widebar{\varepsilon}(g)\cdot e_{i,i+1}(h_i-g_i)) = h
\]\qed

If two 1-cells in $\Theta$ share a vertex via an elementary fold of the form $e_{i,i+1}(q)$, we will call that a vertex of type-$i$. Each 1-cell has $2$-vertices, one of each type.

Let $g\in \Q^2$ so that $g\cdot \mathfrak{C}\subseteq \Theta$ is a 1-cell, and let $\theta\subseteq\{1,2\}$. Then define the following two notations.
\[
g^\theta = (g_i\cdot\chi_{\theta}(i))_i
\]
\[
S(g) = \sum_{\theta\subseteq\{1,2\}} (-1)^{|\theta|} (g^\theta\cdot \mathfrak{C})
\]
That is, $g^\theta$ takes the value of $g$ in any index which is included in $\theta$, and is zero otherwise.
\begin{lemma}
$S(g)\in H_{1}(\Theta; \Q)$.
\end{lemma}
\Proof To prove this, we will pair up 1-cells in $\Theta$ for each $i\in\{1,2\}$ such that the pairs share a vertex of type $i$ and have opposite orientation. In this way, we will see that $S(g)$ has no boundary component and hence is a 1-cycle.

Fix $i\in\{1,2\}$ and let $\omega$ be the other element. Then $g^{\omega\cup\{i\}}$ and $g^\omega$ differ only in the $i$th position, so they share a vertex of type $i$. Furthermore, $|\omega| = |\omega\cup\{i\}|-1$, so they have opposite orientation in the summation of $S(g)$. We can make this argument for each $i$, giving us pairings over each vertex type.\qed

\begin{lemma}
\label{lem:genSet}
$\{S(g)\ |\ g\in \Q^{2}\}$ is a generating set for $H_{1}(\Theta; \Q)$.
\end{lemma}
\Proof Let $R\in H_{1}(\Theta; \Q)$, then we can write
\[
R = \sum_{j\in J}\alpha_j\cdot (g_j \mathfrak{C})
\]
where $J$ is a finite index set, and each $0\neq\alpha_j\in \Q$ and $g_j\in \Q^{2}$. We will show that
\[
R = \sum_{j\in J} \alpha_j\cdot S(g_j)
\]

For $\theta\subseteq\{1,2\}$ and $r\in J$, let $F(\theta, r) = \{j\in J\ |\ g_r^\theta = g_j^\theta\}$.

Let $|\theta|=1$ and fix some $r\in J$. Then the edges $\{g_j\cdot \mathfrak{C}\ |\ j\in F(\theta,r)\}$ all share a common vertex. Notice that if $j\in J$ and $g_j\mathfrak{C}$ also includes this vertex, then $j\in F(\theta, r)$. So, as $R$ is a cycle, and $F(\theta,r)$ indexes over every edge which includes this vertex, it must be the case that $\sum_{j\in F(\theta,r)}\alpha_j = 0$. In addition, this means that there is a subset $J_1^\theta\subseteq J$ such that
\[
J = \bigsqcup_{r\in J^\theta_1} F(\theta,r)
\]

Hence,
\[
\sum_{j\in J}\alpha_j \cdot(g_j^\theta\mathfrak{C}) = \sum_{r\in J_1^\theta} \sum_{j\in F(\theta,r)} \alpha_j\cdot (g_j^\theta\mathfrak{C}) = \sum_{r\in J_1^\theta} \left[\left(\sum_{j\in F(\theta,r)} \alpha_j\right) \cdot (g_r^\theta\mathfrak{C})\right] = \sum_{j\in J_1^\theta} 0 = 0
\]

Notice this also implies that $\sum_{j\in J}\alpha_j = 0$.

That is, we have shown that when $|\theta|$ is 0 or 1, then $\sum_{j\in J} \alpha_j\cdot(g_j^\theta\mathfrak{C}) = 0$. As $g_j^{\{1,2\}} = g_j$ for all $j\in J$. We can combine all of this and see
\begin{eqnarray*}
\sum_{j\in J} \alpha_j\cdot S(g_j) &=& \sum_{j\in J}\alpha_j\sum_{\theta\subseteq\{1,2\}} (-1)^{|\theta|}(g_j^\theta\mathfrak{C})\\
&=& \sum_{\theta\subseteq\{1,2\}} (-1)^{|\theta|} \sum_{j\in J} \alpha_j (g_j^\theta\mathfrak{C})\\
&=& \left(\sum_{\theta\subsetneq\{1,2\}} (-1)^{|\theta|} \sum_{j\in J} \alpha_j (g_j^\theta\mathfrak{C})\right) + \left(\sum_{\theta=\{1,2\}} (-1)^{|\theta|} \sum_{j\in J} \alpha_j (g_j^\theta\mathfrak{C})\right)\\
&=& 0 + \sum_{j\in J} \alpha_j (g_j\mathfrak{C}) = R
\end{eqnarray*}
\qed

Now we define the following functional from edges in $\Theta$ to $\Q$:
\[
\varphi(h\cdot \mathfrak{C}) = h_1h_2
\]
where we are taking $h\in \Q^{2}$ and $h_i$ is the $i$th coordinate of $h$. We can extend this linearly to be a $1$-cochain on the $1$-chains in $\Theta$. As $\Theta$ is $1$-dimensional, $\varphi$ must in fact represent a class in $H^{1}(\Theta; \Q)$.

\begin{lemma}
\label{lem:phiUinvariant}
$\varphi$ is a $U-invariant$ cocycle. That is, for any $u\in U$, and $R\in H_1(\Theta; \Q)$, $\varphi(uR) = \varphi(R)$.
\end{lemma}
\Proof Let $g\in\Q^2$, then $S(g) = \sum_{\theta\subseteq\{1,2\}} (-1)^{|\theta|} (g^\theta\cdot \mathfrak{C})$ and
\[
\varphi(S(g)) = \sum_{\theta\subseteq\{1,2\}} (-1)^{|\theta|}(g^\theta)_1(g^\theta)_2 = g_1g_2
\]
As when $\theta\neq\{1,2\}$, the terms are zero. Let $\varepsilon(u) = (u_1,u_2)\in\Q^2$, then $u\cdot (h\cdot\mathfrak{C}) = ((u_1,u_2) + h)\cdot\mathfrak{C}$ for any $h\in\Q^2$. Now we calculate
\begin{eqnarray*}
\varphi(u\cdot S(g)) &=& \varphi\left(\sum_{\theta\subseteq\{1,2\}} (-1)^{|\theta|} u(g^\theta\cdot \mathfrak{C})\right)\\
&=& \sum_{\theta\subseteq\{1,2\}} (-1)^{|\theta|} (u_1+(g^\theta)_1)(u_2+(g^\theta)_2)\\
&=& u_1u_2 - (u_1+g_1)u_2 - u_1(u_2+g_2) + (u_1+g_1)(u_2+g_2) = g_1g_2
\end{eqnarray*}
So, $\varphi(uS(g)) = \varphi(S(g))$ for any $u\in U$ and $g\in\Q^2$. Lemma~\ref{lem:genSet} tells us that the $S(g)$ are a generating set for the cycles $H_1(\Theta; \Q)$. And if the action of $u$ is invariant on a generating set, then it is invariant for any cycle $R$.\qed\\

\subsection{A local cycle which evaluates to be nonzero.}
\label{sec:localCycle}
In this section, we detail a 1-cycle in $\Xi$ which evaluates to be nonzero under $\varphi\circ\pi$.

Let $P<\SL_3(\Q)$ be the parabolic subgroup which fixes the subspace $\langle e_1,e_2\rangle$ that we defined in \SecRef{sequenceOfVertices}. And again, let $\Omega\subseteq\Xi$ be the open hemisphere complex opposite $P$. We have already shown that this open hemisphere complex consists of the chambers $P\cdot\{D\subseteq\A\ |\ \langle e_3\rangle$ is a vertex of $D\}$.

We can extend the map $\pi\colon U\cdot \mathfrak{C}'\to\Theta$ to the map
\[
\pi\colon \Omega\to [U,U]\setminus\Omega
\]
with the ultimate goal of constructing a cocycle on $\Omega$ by applying $\pi$ to cycles in $\Omega$ and then evaluating $\varphi$ on the image. For this to work, it must be the case that for any cycle $H\in H_1(\Omega)$, $\pi(H) \subseteq \Theta$. We will prove this in \SecRef{reducingDisks}. In this subsection, we construct a single cycle $\mathfrak{B}\in H_1(\Omega; \Q)$ with the property that $\pi(\mathfrak{B})\subseteq \Theta$, and that $\varphi(\pi(\mathfrak{B}))\neq 0$.

Recall that we can realize the chambers of $\Xi$ as complete flags in $\Q^3$, and that any linearly independent basis for $\Q^3$ gives an apartment in $\Xi$. With this in mind, consider the linearly independent vectors:
\[
v_1 = \begin{pmatrix} 0\\0\\1\end{pmatrix}, v_2 = \begin{pmatrix} 0\\1\\1\end{pmatrix}, v_3 = \begin{pmatrix} 1\\1\\1\end{pmatrix}
\]

These determine an apartment in $\Xi$, which we will call $\mathfrak{B}$.

Consider the following two chambers in $\Xi$, which we give by their corresponding complete flags:
\[
C_0 = \mathfrak{C}' \text{ given by } \langle e_3\rangle < \langle e_2,e_3\rangle
\]
\[
C_1 \text{ given by } \langle e_3\rangle < \langle e_1,e_3\rangle
\]

Notice that $C_0,C_1 < \Omega$. Also let:
\[
f_1 = \begin{psmallmatrix} 1 & 1 & 0\\ 0 & 1 & 0\\ 0 & 0 & 1\end{psmallmatrix}\text{ and } f_2 = \begin{psmallmatrix} 1 & 0 & 0\\ 0 & 1 & 1\\ 0 & 0 & 1\end{psmallmatrix}
\]

Using these chambers and elementary matrices, we construct the apartment $\mathfrak{B}$ as a 1-cycle below and seen in \figref{cycleDiagram}.
\[
\mathfrak{B} = C_0 - f_1C_0 + f_1f_2C_0-f_1f_2f_1^{-1}C_1 + f_2C_1 - f_2C_0
\]
\indent We can apply the projection $\pi$.
\[
\pi(\mathfrak{B}) = [C_0] - [f_1C_0] + [f_1f_2C_0] - [f_2C_0]
\]
which is a 1-cycle in $\Theta$. Evaluating by $\varphi$ yields
\[
\varphi(\pi(\mathfrak{B})) = \varphi([C_0]) - \varphi([f_1C_0]) + \varphi([f_1f_2C_0]) - \varphi([f_2C_0])=1
\]

\begin{figure}[t]
	\centerline{\includegraphics{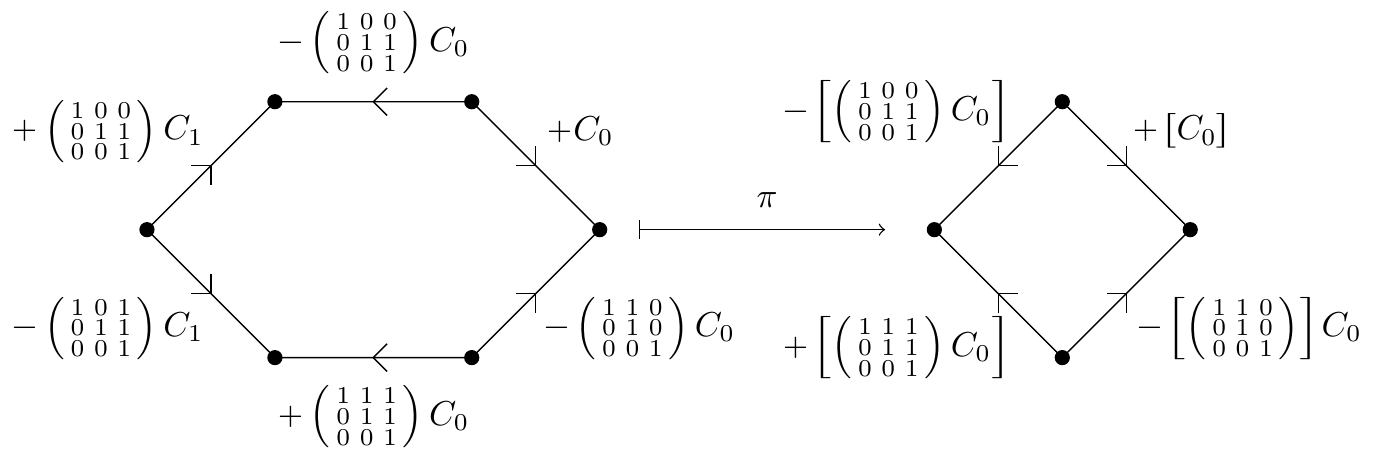}}
	\caption[A nonzero local cycle]{The cycle $\mathfrak{B}$ on the left collapses to a cycle with four edges after applying $\pi$.}
	\figlabel{cycleDiagram}
\end{figure}

\section{A sequence of subgroups of $U(\SL_3(\Q[t]))$}
\label{sec:congruenceSubgroups}
In this \sectionChapter{}, we construct a subgroup $\Congk$ for each $k>0$ that we will associate to $y_k$. 

Consider the map $\epsilon\colon U(\SL_3(\Q[t]))\to (\Q[t])^{2}$ exactly as we did in the local cocycle section, with the kernel being the commutator subgroup of $U(\SL_3(\Q[t]))$. Let $D_k = D_{y_k}$ be the diagonal element we defined in \SecRef{preliminaries}. That is, $y_k = D_k\cdot x_0$. Let's write $D_k = D(t^{a_1(k)},t^{a_{2}(k)},1)$, where $a_1(k)> a_2(k)>0$ for each $k>0$. Now let:
\[
\Congk= \epsilon^{-1}(t^{a_1(k)-a_2(k)+1}\Q[t] \times t^{a_{2}(k)+1}\Q[t])
\]
Let $\pi_k\colon X\to \Congk\setminus X$ be projection. As $\Congk$ contains the commutator subgroup of $U(\SL_3(\Q[t]))$, it is normal in $U(\SL_3(\Q[t]))$. So $U(\SL_3(\Q[t]))$ acts on $\pi_k(X)$, and for any $u\in U(\SL_3(\Q[t]))$ and $x\in X$
\[
\pi_k(u\cdot x) = u\cdot\pi_k(x)
\]

\section{Reducing Disks in $X$}
\label{sec:reducingDisks}

Let $k>0$. In this \sectionChapter{}, we show that for any disk $D\subseteq X$ such that $\partial D \subseteq X_{0}$, the intersection $\pi_k(D)\cap \St(\pi_k(y_k))$ is contained in the image under $\pi_k$ of the translates of the downward chamber by the stabilizer of $y_k$ in $\SL_3(\Q[t])$.

Recall that the elementary matrices in $\{e_{i,j}(t^r)\ |\ j>i\ \text{and}\ r>0\}\subseteq U(\SL_3(\Q[t]))$ are in one-to-one correspondence with the walls of $\A$ which intersect with the interior of $\Sect$. Let $D,D'<\Sect$ be two chambers which share a face $F$. Then there exists some $i,j$ and $r$ such that $F$ is contained in the wall corresponding to $e_{i,j}(t^r)$. The additive subgroup $\{e_{i,j}(qt^r)\ |\ q\in\Q\}\cong\Q$ acts on the chambers of $X$ which contain $F$ as a wall, and the chambers $D$ and $D'$ are a strict fundamental domain for this action.

Now let
\[
d_1=D(t,1,1)
\]
\[
d_2 = D(t,t,1)
\]
Then the convex hull of $\{d_m^l\cdot y_k\}_{l\ge0}$ for $m\in\{1,2\}$ and $k>0$ gives a ray starting at $y_k$ and going off towards a vertex of $\Sect(\infty)\subseteq\A(\infty)$. Let $W$ be a wall corresponding to the elementary matrix $e_{i,j}(t^r)$ such that $y_k\in W$ and the ray is not parallel to $W$. Then $d_m^l\cdot y_k$ is contained in the wall given by $e_{i,j}(t^{r+l})$. This process gives us every wall which the ray passes through.

Now let $\Sect_k\subseteq\A$ be the sector starting at $y_k$ and going off towards $\Sect(\infty)$. Let $\{W^1_k,W^2_k\}$ be the walls of $\Sect_k$. Each gives a half-space $H^i_k\subseteq\A$ which includes $x_0$. That is, the intersection of the complements of the $H_k^i$ gives the interior of $\Sect_k$. Notice that $\left(H^1_k\cap H_k^2\right)\cap\St(y_k)$ gives a single chamber of $\St(y_k)$, let's call this chamber $C_k$. Let $E\subseteq (\St(y_k)\cap\A)$ be a chamber not equal to $C_k$, then fix $i$ such that $E\nsubseteq H_k^i$. Notice that the interior of $H_k^i$ contains exactly one vertex of $C_k$, let's call this vertex $v_i$. Then the ray starting at $v_i$ going through $y_k$ gives a vertex in $\A(\infty)$ which we will call $V_i$. Let $R$ be a ray starting in the interior $E$ and going off towards $V_i$. Notice that $R$ is contained entirely in the complement of $H^i_k$, and is parallel to, but not contained in, $W_k^j$ for $j\neq i$. \figref{reduceDisk} gives a schematic of this arrangement.

\begin{lemma}
\label{lem:wallsExist}
Let $W$ be a wall which $R$ passes through. Then there exists some $r,i,j$ such that:
\begin{enumerate}[i)]
\item $W$ is the wall corresponding to the fold $e_{i,j}(t^r)$, and
\item $\{e_{i,j}(qt^r)\ |\ q\in\Q\}\subseteq\Lambda_k$
\end{enumerate}
\end{lemma}

\Proof First we remark that the walls which $R$ passes through is contained in the set of walls which the parallel ray starting at $y_k$ passes through. Also note that $\Lambda_k$ includes $e_{i,j}(\Q\cdot t^r)$ for each subgroup which corresponds to a wall passing through $\Sect_k$ except for $W_k^1$ and $W_k^2$. We chose $R$ so that it passes only through walls which intersect $\Sect_k$ and does not pass through either $W_k^1$ or $W_k^2$.\qed

\begin{figure}[t]
	\centerline{\includegraphics{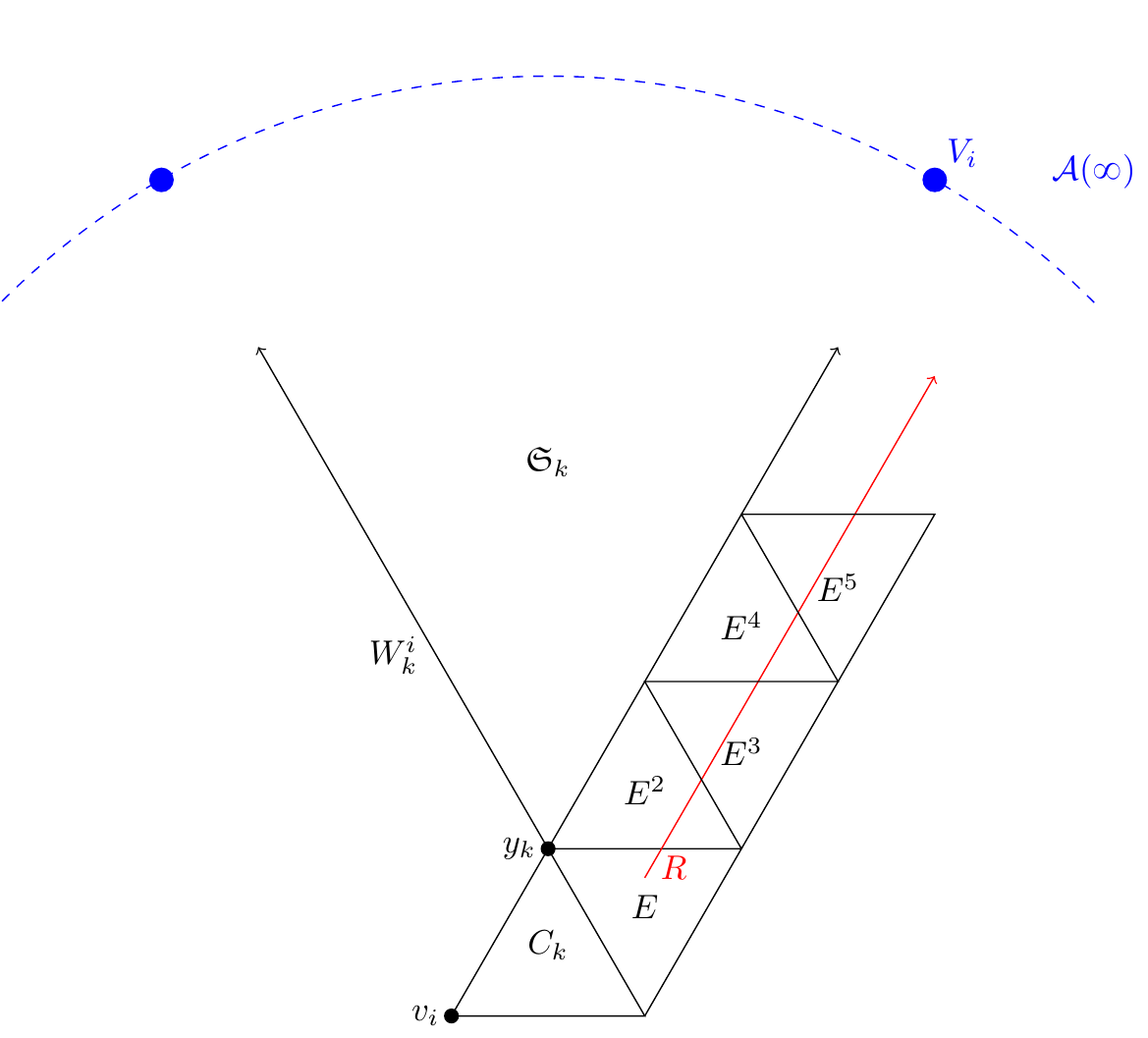}}
	\caption[Reducing disks in $X$]{$R$ passes through an infinite number of chambers where each consecutive pair is a fundamental domain for the star about their shared wall.}
	\figlabel{reduceDisk}
\end{figure}

\begin{lemma}
Let $k>0$ and $D\subseteq X$ be a 2-chain such that $\partial D\subseteq X_0$. Then $\pi_k(D)$ is a 2-chain in $\pi_k(X)$ and $\partial \pi_k(D) \cap \pi_k(\Sect_k) = \varnothing$. 
\end{lemma}
\Proof Recall that we chose $y_k$ in such a way that $\St(y_k)\cap \C = \varnothing$ and $\Sect_k\subseteq\Sect$. Hence $\Sect_k\cap \C = \varnothing$ and furthermore, as $X_0\cap\Sect = \C$ , it is the case that $\Sect_k\cap X_0 = \varnothing$. As $\Sect$ is a strict fundamental domain for the action of $\SL_3(\Q[t])$ on $X$, and $\Congk < \SL_3(\Q[t])$, we can conclude that $\pi_k(\Sect_k)\cap\pi_k(X_0) = \varnothing$. As $\partial D\subseteq X_0$, we have $\partial\pi_k(D)\subseteq\pi_k(X_0)$, which is disjoint from $\pi_k(\Sect_k)$.\qed

\begin{lemma}
Let $D\subseteq X$ be a 2-chain such that $\partial D\subseteq X_0$, then the support of $\pi_k(D)$ in $\St(\pi_k(y_k))$ is contained in $\pi_k(U(\SL_3(\Q[t]))\cdot C_k)$.
\end{lemma}
\Proof We chose the $y_k$ in such a way that $\St(y_k)\cap (\SL_3(\Q[t])\cdot\C) = \varnothing$ for each $k$. Furthermore, $\Sect_k\cap(\SL_3(\Q[t])\cdot\C) = \varnothing$ for each $k$. Note that $X_0\subseteq (\SL_3(\Q[t])\cdot\C)$, but using this larger space for possible boundary components will make the proof easier.

Let $E$ be a chamber of $\St(y_k)\cap\A$ which is not equal to $C_k$. We can construct a ray $R$ as above starting in the interior of $E$. Let $E=E^1,E^2,\dots$ be the sequence of chambers which $R$ passes through, so that $E^i$ and $E^{i+1}$ share a face. By lemma~\ref{lem:wallsExist}, we know that $\pi_k(E^i)$ and $\pi_k(E^{i+1})$ are the only two chambers which include that common face in $\pi_k(X)$. As the ray $R$ stays within a half-space that does not include any chambers of $(\SL_3(\Q[t])\cdot\C)$, we know that it never passes through a boundary component of $D$. So, if $\pi_k(E)$ is in the support of $\pi_k(D)$, then in fact, every $\pi_k(E^i)$ must also be in the support. Each maps to a unique chamber in $\pi_k(X)$, as each is a distinct chamber in $\Sect$, a fundamental domain for the action of $\SL_3(\Q[t])$ on $X$. This contradicts the fact that $D$ is compact.

Now suppose that $E'$ is any chamber in $\St(y_k)$ that is not in the $\SL_3(\Q[t])$ orbit of $C_k$. Then $E' = u\cdot E$ for some $u\in \SL_3(\Q[t])$ and $E\in\St(y_k)\cap\A$. Then $u^{-1}D$ is a disk with boundary in $\SL_3(\Q[t])\cdot\C$, and the above paragraph tells us that $\pi_k(u^{-1}E')$ cannot be contained in the support of $\pi_k(u^{-1}D)$. Hence, $\pi_k(E')$ cannot be in the support of $\pi_k(D)$.\qed\\

\section{A Family of Global Cocycles}
\label{sec:globalCocycles}

In this \sectionChapter{}, we construct an infinite family of independent cocycles in $H^2(\Gamma\setminus Y; \Q)$.

\subsection{A Function on Global disks}
In this section, we define a function $\varphi_k$ for each $k>0$ that takes 2-disks in $Y$ to $\Q$.

In \SecRef{localStars}, we defined $\eta_k$, which took chambers in $\St(y_k)$ to $\Xi$ and an equivariant map $\rho_k\colon U(\Stab_{\SL_3(\Q[t])}(y_k))\to U(\SL_3\Q)$.

\begin{lemma}
For any $k>0$, $\rho_k((\Congk)_{y_k}) = [U(\SL_3\Q), U(\SL_3\Q)]$.
\end{lemma}
\Proof This follows directly from the construction of $\Congk$ and lemma~\ref{lem:commutators}.\qed\\

In constructing the map $\eta_k$, we chose to send $C_k\subseteq\St(y_k)$ to $\mathfrak{C}'\subseteq\Xi$. This gives us the connection between $\Congk$ and $\Theta$ in the next lemma. 

\begin{lemma}
There is an isomorphism
\[
\widebar{\eta}_k\colon \pi_k(\Stab_{\SL_3(\Q[t])}(y_k)\cdot C_k) \to \pi(U(\SL_3\Q)\cdot\mathfrak{C}')=\Theta
\]
such that for any $u\in \Stab_{\SL_3(\Q[t])}(y_k)$ and cell $E \subset \pi_k(\Stab_{\SL_3(\Q[t])}(y_k)\cdot C_k)$ we have:
\[
\rho_k(u)\cdot\widebar{\eta}_k(E) = \widebar{\eta}_k(u\cdot E)
\]
\end{lemma}
\Proof This map is induced by $\eta_k$. Taking the quotient $\pi_k$ in the domain is equivalent to taking the quotient $\pi$ in the image, as seen by lemma 9.1 above. Thus this induced map is still an isomorphism, and can be written as:
\[
\widebar{\eta}_k(\pi_k(u\cdot C_k)) = \pi(\eta_k(u\cdot C_k))
\]
where $u\in\Stab_{\SL_3(\Q[t])}(y_k)$. Using the fact that $\rho_k$ and $\eta_k$ are equivariant, we can manipulate this:
\begin{eqnarray*}
\rho_k(u)\cdot\widebar{\eta}_k(E) &=& \rho_k(u)\cdot\pi(\eta_k(E))\\
&=& \pi( \rho_k(u)\cdot\eta_k(E))\\
&=& \pi(\eta_k(u\cdot E))
\end{eqnarray*}\qed

Given a disk $D$ with boundary in $X_0$, we showed above that for any $k>0$, $\pi_k(D)\cap\St(\pi_k(y_k))$ is contained in $\pi_k(\Stab_{\SL_3(\Q[t])}(y_k)\cdot C_k)$. So we can construct a map $\Eta_k$ that takes such disks to $\Theta$, which is the composition of $\pi_k$, restriction to $\St(\pi_k(y_k))$ and $\widebar{\eta}_k$. We remark that this map $\Eta_k$ is dimension reducing, as it includes the map $\downarrow_k$ which takes chambers in $\St(y_k)$ to their boundaries in $\Lk(y_k)$.

\begin{lemma}
Let $D$ be a $2$-disk with $\partial D\subseteq X_{0}$ and $k>0$, then $\Eta_k(D)$ is a $1$-cycle in $\Theta$.
\end{lemma}
\Proof This follows from the fact that $\pi_k(D)$ has no boundary components in $\St(\pi_k(y_k))$.\qed\\

Now we proceed to construct the map
\[
\varphi_k\colon \{\text{2-disks } D \subseteq X\ |\ \partial D\subseteq X_0\} \to \Q
\]
for each $k>0$. 

From construction, we have two types of 2-cells in $Y$: those that are contained in $X_{0}$ and those that we glued in to construct $Y$ from $X_{0}$. For cells of the first type, we will have $\varphi_k$ take them to 0. For $E$, a cell of the second type, define:
\[
\varphi_k(E) = \varphi(\Eta_k(\psi(E)))
\]
where $\varphi$ is the local cocycle defined in \SecRef{localCocycles}.

\begin{lemma}
\label{lem:invariance}
Let $E < Y$ be a 2-cell which is not contained in $X_0$. Then for any $u\in U(\Gamma)$, $\varphi_k(uE) = \varphi_k(E)$.
\end{lemma}
\Proof Any $u\in U(\Gamma)$ can be decomposed into $u=u_0u_1$ where $u_0\in\Congk$ and $u_1\in\Gamma_{y_k}$. Hence
\[
\varphi_k(uE) = \varphi(\Eta_k(\psi(u_0u_1E))) = \varphi(\Eta_k(u_0u_1\psi(E)))
\]
as $\psi$ is $\Gamma$-equivariant. And as $u_0\in\Congk$ and $u_1\in \Gamma_{y_k}<\Stab_{\SL_3(\Q[t])}(y_k)$ which acts on $\St(\pi_k(y_k))$, we have:
\[
\varphi(\Eta_k(u_0u_1\psi(E))) = \varphi(\rho_k(u_1)\cdot\Eta_k(\psi(E)))
\]
$\rho_k(u_1)\in U(\SL_3(\Q))$, and $\varphi$ is $U(\SL_3(\Q))$-invariant. So in fact, we have:
\[
\varphi(\rho_k(u_1)\cdot\Eta_k(\psi(E))) = \varphi(\Eta_k(\psi(E))) = \varphi_k(E)
\]
\qed

\begin{lemma}
\label{lem:wellDefinedSum}
Let $E$ be a 2-cell of the second type in $Y$, then the summation
\[
\sum_{\gamma U(\Gamma)\in\Gamma / U(\Gamma)}\varphi_k(\gamma^{-1} E)
\]
is well-defined.
\end{lemma}

\Proof Lemma~\ref{lem:invariance} shows that each term in the summation is invariant with respect to choice of coset representative. So all that needs to be shown is that only finitely many terms in the sum are non-zero.
Suppose that $\gamma,\gamma'\in\Gamma$ and $\pi_k(\gamma y_k) = \pi_k(\gamma' y_k)$. Suppose that $\gamma,\gamma'\in\Gamma$ and $\pi_k(\gamma y_k) = \pi_k(\gamma' y_k)$. Then $\gamma y_k = \lambda \gamma' y_k$ for some $\lambda\in\Lambda_k$. That is, $\gamma^{-1}\lambda\gamma'=u$ for some $u\in\Stab_{\SL_3(\Q[t])}(y_k)$. Now we manipulate this below, using the fact that $\lambda\in B(\SL_3(\Q[t])$, a normal subgroup of $\SL_3(\Q[t])$.
\begin{eqnarray*}
\gamma^{-1}\lambda\gamma' &=& u\\
\gamma^{-1}\lambda\gamma'(\gamma'^{-1}\gamma)\gamma^{-1}\gamma' &=& u\\
\gamma^{-1}\lambda \gamma (\gamma^{-1}\gamma') &=& u\\
\lambda' \gamma^{-1}\gamma' &=& u\\
\gamma^{-1}\gamma' &=& u\lambda'^{-1}
\end{eqnarray*} 
where $\lambda^{-1}\in B(\SL_3(\Q[t])$. $u$ and $\lambda'^{-1}$ are both upper-triangular, hence so is their product. So $\gamma^{-1}\gamma'\in(\Gamma\cap B(\SL_3(\Q[t]))) = U(\Gamma)$.

Thus, if $\gamma,\gamma'$ represent distinct cosets, $\gamma\St(\pi_k(y_k))\cap\gamma'\St(\pi_k(y_k)) = \varnothing$, as we are talking about the open star. For $\varphi_k(E)$ to evaluate to be non-zero, $\pi_k(E)$ must include a chamber in $\St(\pi_k(y_k))$. Thus if infinitely many terms were non-zero, $E$ would include infinitely many distinct chambers in its support (as each star is disjoint) contradicting the fact that $E$ is a compact disk.\qed\\

\subsection{A $2$-cochain on $\Gamma\setminus Y$}
Let $\Gamma E$ be a $2$-cell in $\Gamma\setminus Y$. For each $k>0$, we define the following function to $\Q$
\[
\Phi_k(\Gamma E) = \sum_{\gamma U(\Gamma)\in\Gamma / U(\Gamma)}\varphi_k(\gamma^{-1} E)
\]
which we can extend linearly to be a $2$-cochain on $\Gamma\setminus Y$.
\begin{lemma}
$\Phi_k$ is well-defined and represents a class in $H^{2}(\Gamma\setminus Y; \Q)$.
\end{lemma}
\Proof We showed in lemma~\ref{lem:wellDefinedSum} that this sum is finite for any cell $E$, so to show that $\Phi_k$ is well-defined, we have to show that it is invariant with respect to the coset representative of $\Gamma E$. That is, we have to show that
\[
\Phi_k(\Gamma E) = \Phi_k(\Gamma h E)
\]
for any $2$-cell $E$ in $Y$ and $h\in\Gamma$. We prove this directly.
\begin{eqnarray*}
\Phi_k(\Gamma h E) &=& \sum_{\gamma U(\Gamma)\in\Gamma / U(\Gamma)}\varphi_k(\gamma^{-1} hE)\\
&=& \sum_{h\gamma U(\Gamma)\in\Gamma / U(\Gamma)}\varphi_k((\gamma^{-1}h^{-1}) hE)\\
&=& \sum_{h\gamma U(\Gamma)\in\Gamma / U(\Gamma)}\varphi_k(\gamma^{-1}E)\\
&=& \sum_{\gamma U(\Gamma)\in\Gamma / U(\Gamma)}\varphi_k(\gamma^{-1}E) = \Phi_k(\Gamma E)
\end{eqnarray*}
where we can make these substitutions because we showed in lemma~\ref{lem:invariance} that the value of each term is invariant up to choice of coset representative.

To show that $\Phi_k$ represents a class in cohomology, we must show that it is in the kernel of the coboundary map. This is equivalent to showing that it evaluates the boundary of $3$-chains in $\Gamma\setminus Y$ to be zero.

Any $3$-cell in $\Gamma\setminus Y$ is the image of a $3$-cell in $Y$. Let $D^3$ be such a $3$-cell in $Y$. $\psi$ takes the boundary, a $2$-sphere into $X$, a $2$-dimensional contractible space. Thus, as a cycle, $\psi(\partial D_3) = 0$, and therefore, $\Phi_k(\partial(\Gamma D_3)) = 0$.\qed\\

\begin{remark}
\label{rem:rebuttal}
We are now in a position to examine why the similarly defined cocycles in \cite{CK} each evaluate to be zero. Suppose that we chose $\Gamma$ to be all of $\SL_3(\Z[t])$. If there exists a $\gamma U(\Gamma) \in \Gamma\setminus U(\Gamma)$ such that $\varphi_k(\gamma^{-1}E)\neq 0$ for some 2-cell $E$, then $D(-1,-1,1)\gamma$, $D(1,-1,-1)\gamma$, and $D(-1,1,-1)\gamma$ are also such elements for $H$, as these diagonal subgroups each fix $y_k$. But, we notice that $\varphi_k(\gamma^{-1}E) = \varphi_k(\gamma^{-1} D(-1,-1,1)E) = -\varphi_k(\gamma^{-1} D(1,-1,-1)E)=\varphi_k(\gamma^{-1} D(-1,1,-1)E)$. Hence, summing over these four cosets will give zero, and furthermore, we can break up the entire summation from lemma~\ref{lem:wellDefinedSum} into these non-intersecting groups of 4, each of which sums to zero. They do not provide an example of a cycle which is evaluated to be nonzero in \cite{CK} because such a cycle cannot exist.
\end{remark}

\subsection{A Sequence of Cycles}
In this section, we prove the following two things:

\begin{proposition}
\label{prop:cocyclesExist}
For each $k>0$, there exists an $\mathfrak{H}_k\in H_{2}(\Gamma\setminus Y; \Q)$ such that $\Phi_k(\mathfrak{H}_k)\neq 0$ and $\Phi_j(\mathfrak{H}_k)=0$ for all $j>k$.
\end{proposition}

\begin{proposition}
\label{prop:infiniteDimensional}
$H_2(\Gamma\setminus Y; \Q)$ is infinite dimensional.
\end{proposition}

\subsection{A Free Subgroup of $\SL_2(\Z[t])$}
First we take a detour and examine a particular subgroup of $\SL_2(\Z[t])$ isomorphic to $\Z$, which we will later embed in $\SL_3(\Z[t])$.

Let
\[
a = \begin{pmatrix} 1&t\\ 0 & 1 \end{pmatrix}\begin{pmatrix} 1&0\\ t & 1 \end{pmatrix} = \begin{pmatrix} t^2+1&t\\ t & 1 \end{pmatrix}
\]
$a$ has eigenvalues $\frac{1}{2}(2+t^2\pm t\sqrt{t^2+4})$, where $\sqrt{t^2+4}$ represents the Laurent series expansion of $\sqrt{t^2+4}$, an element with valuation -1. These are distinct, hence $a$ is diagonalizable in $\SL_2(\Q((t^{-1})))$. That is, there exists $g\in\SL_2(\Q((t^{-1})))$ such that $gag^{-1}$ is a diagonal matrix. Furthermore, each eigenvalue has nontrivial norm, so $a$ has infinite order and $\langle a\rangle=\Z$. Also, $\langle gag^{-1}\rangle\cong\Z$ is a diagonal subgroup of $\SL_2(\Q((t^{-1})))$. In fact, we can specify a $g$ precisely:
\[
g = \begin{pmatrix}
\frac{1}{2}(t-\sqrt{t^2+4}) & \frac{-t-\sqrt{t^2+4}}{2\sqrt{t^2+4}}\\
1 & -\frac{1}{\sqrt{t^2+4}}
\end{pmatrix}
\]
and can confirm from observation and expansion that $g\in\SL_2(\Q[[t^{-1}]])$. 

Elements of $\langle a\rangle$ act on vectors in $(\Z[t])^2$ in the standard fashion.
\begin{lemma}
$\langle \{a^k\cdot\Z^2\}_{k\in\Z}\rangle = (\Z[t])^2$.
\end{lemma}
\Proof $(\Z[t])^2$ is additive, so it is sufficient to show that $\Vect{t^r}{0}, \Vect{0}{t^r}\in\langle \{a^k\cdot\Z^2\}_{k\in\Z}\rangle$ for all $r\ge0$. We do this by induction. $a^0\Vect{1}{0} = \Vect{1}{0}$ and $a^0\Vect{0}{1}=\Vect{0}{1}$, so $\Vect{1}{0},\Vect{0}{1}\in\langle \{a^k\cdot\Z^2\}_{k\in\Z}\rangle$, so we have the base case.

Before proceeding to the inductive step, notice that for $\ell>0$:
\[
a^\ell = \begin{pmatrix} t^{2\ell} + a_{11}(t) & t^{2\ell-1}+a_{12}(t)\\ t^{2\ell-1} + a_{21}(t) & t^{2\ell-2}+a_{22}(t)\end{pmatrix}
\]
\[
a^{-\ell} = \begin{pmatrix} t^{2\ell-2} + b_{11}(t) & -t^{2\ell-1}+b_{12}(t)\\ -t^{2\ell-1} + b_{21}(t) & t^{2\ell}+b_{22}(t)\end{pmatrix}
\]
where:
\begin{eqnarray*}
\deg(a_{11}),\deg(b_{22}) &<& 2\ell\\
\deg(a_{12}),\deg(a_{21}),\deg(b_{12}),\deg(b_{21}) &<& 2\ell-1\\
\deg(a_{22}),\deg(b_{11}) &<& 2\ell-2\\
\end{eqnarray*}
So we will in fact prove two steps of induction at a time. That is, let $r>0$ and assume that the $\Vect{t^s}{0},\Vect{0}{t^s}\in W$ for $s<2r-1$. We will show that it is also the case when $s=2r-1,2r$.\\
Consider:
\[
a^r\cdot\Vect{0}{1} = \Vect{t^{2r - 1}+a_{12}(t)}{t^{2r-2}+a_{22}(t)}
\]
By our inductive assumption, we know that $\Vect{a_{12}(t)}{t^{2r-2}+a_{22}(t)}\in \langle \{a^k\cdot\Z^2\}_{k\in\Z}\rangle$ so we have that $\Vect{t^{2r-1}}{0}\in \langle \{a^k\cdot\Z^2\}_{k\in\Z}\rangle$. Observing $a^{-r}\cdot \Vect{-1}{0}$ gives us that $\Vect{0}{t^{2r-1}}\in W$ as well.\\
Now consider:
\[
a^r\cdot\Vect{1}{0} = \Vect{t^{2r}+a_{11}(t)}{t^{2r-1}+a_{21}(t)}
\]
and by the same logic, we find that $\Vect{t^{2r}}{0}\in W$. Then looking at $a^{-r}\cdot \Vect{0}{1}$ shows that $\Vect{0}{t^{2r}}\in \langle \{a^k\cdot\Z^2\}_{k\in\Z}\rangle$ as well. This completes the proof.\qed

\begin{lemma} No nontrivial element of $\langle a\rangle$ fixes a point in the Euclidean building for $\SL_2(\Q((t^{^-1})))$.
\end{lemma}
\Proof We saw above that $a$ is diagonalizable and hence eigenvalues with nontrivial norm. Hence $a$ has a hyperbolic translation of the building. Such translations have no fixed points. This reasoning applies to every nontrivial power of $a$.\qed\\

We summarize our findings in the lemma below.

\begin{lemma}
For the element $a\in\SL_2(\Z[t])$ described above, the following hold:
\begin{enumerate}[i)]
\item $\langle a\rangle \cong\Z$
\item There is some $g\in\SL_2(\Q[[t^{-1}]])$ such that $gag^{-1}$ is diagonal
\item No nontrivial element of $\langle a\rangle$ fixes a point in the Euclidean building for $\SL_2(\Q((t^{-1})))$
\end{enumerate}
\end{lemma}

\subsection{A Small Unipotent Toolkit}
In this subsection, we develop a small unipotent toolkit that will be used later in the paper. We borrow most relevant terms from \cite{BoSe}.

Let $\Delta=\{\alpha_{12},\alpha_{23}\}$ be the set of simple roots of the diagonal maximal torus\\
in $\SL_3(\Q((t^{-1})))$. Each $\alpha\in\Delta$ gives several subgroups of $\SL_3(\Q((t^{-1})))$:
\begin{itemize}
\item $S_\alpha$, the connected component of the kernel of $\alpha$
\item $Z(S_\alpha)$, the centralizer of $S_\alpha$ in $\SL_3(\Q((t^{-1})))$
\item $U_\alpha$, a unipotent subgroup of $\SL_3(\Q((t^{-1})))$
\item $P_\alpha=Z(S_\alpha)U_\alpha$, a parabolic subgroup of $\SL_3(\Q((t^{-1})))$
\end{itemize}
Let $H$ be any of the groups above and $K<\Q((t^{-1}))$. 
\[
H(K) = H\cap\SL_3(K)
\]
Notice that for either $\alpha\in\Delta$, $Z(S_\alpha)(\Z[t])\cong\GL_2(\Z[t])$, and so we can map elements of $\SL_2(\Z[t])$ to $Z(S_\alpha)(\Z[t])$. Specifically, for $a\in\SL_2(\Z[t])$ as above, let $a_\alpha$ be the image of this element in $Z(S_\alpha)(\Z[t])$. This will be a copy of $a$ in either the upper-left or bottom-right corner of $Z(S_\alpha)$ (depending on $\alpha$). Similarly, let $g_\alpha\in Z(S_\alpha)$ the copy of $g$ from the prior subsection in $Z(S_\alpha)$. Thus we have that $\langle g_\alpha a_\alpha g_\alpha^{-1}\rangle\cong\Z$ is a diagonal subgroup of $\SL_3(\Q((t^{-1})))$.

These subgroups are connected to the Euclidean building $X$ by virtue of the fact that the visual boundary of $X$ (which we denote $X(\infty)$) is the spherical building for $\SL_3(\Q((t^{-1})))$ with simplices corresponding to parabolic subgroups of $\SL_3(\Q((t^{-1})))$. For $\alpha\in\Delta$, the $P_\alpha$ correspond to vertices of the standard chamber of the standard apartment of $X(\infty)$. Vertices in $X(\infty)$ are defined to be an equivalence class of rays in $X$. So for $\alpha\in\Delta$ and a point $p\in X$, let $R_\alpha(p)$ be the ray starting at $p$ which defines $P_\alpha$. That is, it is the ray starting at $p$ and going off towards $P_\alpha$ in the visual boundary. We let $R_\alpha = R_\alpha(x_0)$ for convenience.

We remark in this paragraph on several connections between the $R_\alpha$ and our subgroups above. The convex hull of $S_\alpha\cdot x_0$ is a line which contains $R_\alpha$, and for $h\in P_\alpha$,  $R_\alpha(h\cdot x_0) = h R_\alpha$. Let $\Sect\subseteq\A\subseteq X$ be the standard sector of the standard apartment. Then:
\[
\partial\Sect = R_{\alpha_{12}}\cup R_{\alpha_{23}}\text{   and   } R_{\alpha_{12}}\cap R_{\alpha{23}}=x_0
\]

Let $U(\Z[t])$ be the upper-triangular unipotent subgroup of $\SL_3(\Z[t])$. Then $U_\alpha(\Z[t])<U(\Z[t])$. We have the following:
\begin{lemma}
$U_\alpha(\Z[t])\cdot R_\alpha$ is the connected component of $R_\alpha$ in $U(\Z[t])\cdot R_\alpha$.
\end{lemma}
\Proof First note that each element of $U_\alpha(\Z[t])$ fixes a vertex of $R_\alpha$, so $U_\alpha(\Z[t])\cdot R_\alpha$ is certainly connected. Now we argue that for $u\in U(\Z[t])$, either $u$ fixes no vertices of $R_\alpha$ and hence takes it to a disjoint ray, or there exists some $u_0\in U_\alpha(\Z[t])$ such that $u_0R_\alpha = uR_\alpha$.

Now we have two cases to consider, when 1: $R_\alpha$ is the convex hull of $\{D(t^k,1,1)\cdot x_0\}_{k\ge 0}$ and $U_\alpha(\Z[t]) = \left\{\begin{psmallmatrix} 1 & * & *\\ 0 & 1 & 0\\ 0 & 0 & 1\end{psmallmatrix}\right\}$, and when 2: $R_\alpha$ is the convex hull of $\{D(t^k,t^k,1)\cdot x_0\}_{k\ge 0}$ and $U_\alpha(\Z[t]) = \left\{\begin{psmallmatrix} 1 & 0 & *\\ 0 & 1 & *\\ 0 & 0 & 1\end{psmallmatrix}\right\}$.

We begin with case 1. Suppose that $u=\begin{psmallmatrix} 1 & x & y\\ 0 & 1 & z \\ 0 & 0 & 1\end{psmallmatrix}$ fixes $D(t^k,1,1)x_0$. Let $u_0 = \begin{psmallmatrix} 1 & x & y-xz\\ 0 & 1 & 0 \\ 0 & 0 & 1\end{psmallmatrix}$. Notice that $D(t^{-k},1,1)uD(t^k,1,1)\in\Stab(x_0)$ implies that $z\in\Z$. Now we show that for any $j\ge0$, $D(t^{-j},1,1) u_0^{-1}u_1D(t^j,1,1)\in\Stab(X_0)$, which would imply that $uR_\alpha = u_0 R_\alpha$. We can calculate this, and see that:
\[
D(t^{-j},1,1) u_0^{-1}u_1D(t^j,1,1) = \begin{pmatrix} 1 & 0 & 0 \\ 0 & 1 & z\\ 0 & 0 & 1\end{pmatrix}\in\Stab(x_0)
\]

The same technique works for the second case. If $u=\begin{psmallmatrix} 1 & x & y\\ 0 & 1 & z \\ 0 & 0 & 1\end{psmallmatrix}$, then we let $u_0 = \begin{psmallmatrix} 1 & 0 & -y\\ 0 & 1 & -z \\ 0 & 0 & 1\end{psmallmatrix}$ and can run the same argument.\qed





\begin{lemma}
Let $\alpha\in\Delta$ and nontrivial $u\in U_\alpha$. Then the part of $\A$ which is fixed by $u$ is a cone composed of every ray at a common basepoint which is asymptotically in a chamber in $\A(\infty)$ which has $P_\alpha$ as a vertex. Furthermore, for any point $p\in\A$, all but at most a compact segment of $R_\alpha(p)$ is fixed by $u$. 
\end{lemma}

\subsection{A Subspace of the `Boundary' of the Standard Sector}
In this subsection, we construct a space for each $\alpha\in\Delta$ which is contained in the $\SL_3(\Q[t])$-orbit of the closed star about $\partial\Sect$.

As $a_\alpha$ is diagonalizable, has infinite order, and has no fixed point, the convex hull of $\{a_\alpha^k\cdot x_0\}_{k\in\Z}$ is an infinite geodesic in $X$ which we will call $A_\alpha$. As $a_\alpha$ commutes with $S_\alpha$, we see that together, they span an apartment in $X$, specifically $g_\alpha^{-1}\A$. This apartment contains $A_\alpha$ and the ray $R_\alpha$. Furthermore, by the construction of $a$, we see that a fundamental domain for the action of $\SL_3(\Q[t])$ on $g_\alpha^{-1}\A$ is the closed star about the line in $\A$ which contains $R_\alpha$. This space in turn is in the orbit of the closed star about $\partial\Sect$. To see this, recall that $\Sect$ is a strict fundamental domain for the action of $\SL_3(\Q[t])$ on all of $X$. Now we will show that the space $U_\alpha(\Z[t])\cdot g^{-1}\A$ has some nice properties:

\begin{lemma}
\label{lem:UAalphaConnected}
$U_\alpha(\Z[t])\cdot A_\alpha$ is connected.
\end{lemma}
\Proof The vertices $\{a_\alpha^k\cdot x_0\}_{k\in\Z}$ are contained in $A_\alpha$. If we can show that the $U_\alpha(\Z[t])$-stabilizers of these vertices generate all of $U_\alpha(\Z[t])$, then we have the claim.\\
\indent First note that $U_\alpha(\Z[t])\cong(\Z[t])^2$ by the map which we will call $F$. For each $k\in\Z$, the set $a_\alpha^kU_\alpha(\Z)a_\alpha^{-k}$ is contained in the $U_\alpha(\Z[t])$-stabilizer, so if we can show that these sets generate, then we are done. Well, $F(a_\alpha^kU_\alpha(\Z)a_\alpha^{-k}) = a^{\pm k}\cdot \Z^2$. We showed above that these sets generate all of $(\Z[t])^2$, and as $F$ is an isomorphism, we have that their preimages under $F$ must generate all of $U_\alpha(\Z[t])$.\qed

There is a natural continuous map $\hat{\mathcal{F}}_\alpha\colon g_\alpha^{-1}\A\times [0,\infty)\to g_\alpha^{-1}\A$ which moves any point $p$ along $R_\alpha(p)$ at unit-speed. To see this, remember that any apartment in $X$ is a Euclidean plane and each $R_\alpha(p)$ is parallel. We can extend this map to a continuous map on all of $U_\alpha(\Z[t])\cdot g_\alpha^{-1}\A\times [0,\infty)$ by taking:
\[
\mathcal{F}_\alpha(u\cdot p,s) = u\hat{\mathcal{F}}_\alpha(p,s)
\]
where $u\in U_\alpha(\Z[t])$ and $p\in g_\alpha^{-1}\A$. It remains to show that this map is well-defined.

\begin{lemma}
The map $\mathcal{F}_\alpha$ is well-defined.
\end{lemma}
\Proof First we show that $g_\alpha^{-1}\A$ is a strict fundamental domain for the action of $U_\alpha(\Z[t])$ on $U_\alpha(\Z[t])\cdot g_\alpha^{-1}\A$. That is, that if $p,p'\in g_\alpha^{-1}\A$ and $u,u'\in U_\alpha(\Z[t])$ such that $up=u'p'$, then $p=p'$. In this case, $u^{-1}u'p'=p$. As $g_\alpha$ normalizes $U_\alpha$, we have that $g_\alpha u^{-1}u'g_\alpha^{-1}\in U_\alpha$. Furthermore, $g_\alpha p, g_\alpha p'\in\A$ and $g_\alpha u^{-1}u'g_\alpha^{-1}(g_\alpha p') = g_\alpha p$. $\A$ is contained in the fundamental domain for any unipotent acton on $X$, thus $g_\alpha p' = g_\alpha p$, which implies that $p=p'$.\\
Now we assume that $u,u'\in U_\alpha(\Z[t])$ and $p\in g_\alpha^{-1}\A$ such that $up=u'p$. Then we have to show that for all $s\in[0,\infty)$, $\mathcal{F}_\alpha(up,s) = \mathcal{F}_\alpha(u'p,s)$. That is, that $u\hat{\mathcal{F}}_\alpha(p,s) = u'\hat{\mathcal{F}}_\alpha(p,s)$, or equivalently that $u^{-1}u' \hat{\mathcal{F}}_\alpha(p,s) = \hat{\mathcal{F}}_\alpha(p,s)$.

Notice that $\{\hat{\mathcal{F}}_\alpha(p,s)\}_{s\in [0,\infty)}=R_\alpha(p)$. So it is sufficient to show that $u^{-1}u'$ fixes $R_\alpha(p)$.

Again, we translate to the standard apartment of $X$. As $u^{-1}u'$ fixes $p$, then $g_\alpha u^{-1}u' g_\alpha^{-1}$ fixes $g_\alpha p\in\A$. And again, we know that $g_\alpha u^{-1} u' g_\alpha^{-1}\in U_\alpha$. In the lemma above, we showed that such unipotent elements fix a convex set in $\A$ which includes $P_\alpha$ in its visual boundary. Thus, as $g_\alpha u^{-1} u'g_\alpha^{-1}$ fixes $g_\alpha p$, it also fixes the entire ray $R_\alpha(g_\alpha p)$. This ray is exactly $g_\alpha R_\alpha(p)$, so we have:
\begin{eqnarray*}
u^{-1}u'\cdot R_\alpha(p) &=& (g_\alpha^{-1} g_\alpha) u^{-1} u' (g_\alpha^{-1} g_\alpha)\cdot R_\alpha(p)\\
&=& g_\alpha^{-1} ( g_\alpha u^{-1} u' g_\alpha^{-1}) g_\alpha R_\alpha(p)\\
&=& g_\alpha^{-1} g_\alpha R_\alpha(p) = R_\alpha(p)
\end{eqnarray*}
\qed

Eventually, we want to use our map $\mathcal{F}_\alpha$ to construct homotopies which push loops in $U_\alpha(\Z[t])\cdot A_\alpha$ to loops in $g_\alpha^{-1}\A$. To do this, we need the following:

\begin{lemma}
Let $p\in g_\alpha^{-1}\A$ and $u\in U_\alpha(\Z[t])$. Then there exists $S>0$ such that $\mathcal{F}_\alpha(up,s)\in g_\alpha^{-1}\A$ for $s\ge S$.
\end{lemma}
\Proof $\{\mathcal{F}_\alpha(up,s)\}_{s\ge 0}$ is exactly the ray $R_\alpha(up) = uR_\alpha(p)$. As $g_\alpha$ normalizes $U_\alpha$, there is some $u'\in U_\alpha$ such that $g_\alpha u = u'g_\alpha$. So:
\[
g_\alpha R_\alpha(up) = g_\alpha u R_\alpha(p) = u' g_\alpha R_\alpha(p) = u' R_\alpha(g_\alpha p)
\]
$R_\alpha(g_\alpha p)$ is a ray contained in $\A$, and so we know that all but a compact portion is fixed by $u'$. That is, there exists some $S>0$ such that for all $s\ge S$, $\mathcal{F}_\alpha(g_\alpha p, s)$ is fixed by $u'$. That is, $u'\mathcal{F}_\alpha (g_\alpha p,s) = \mathcal{F}_\alpha(g_\alpha p, s)$. Now we can apply $g_\alpha^{-1}$ to both sides and find that for all $s\ge S$:
\begin{eqnarray*}
g_\alpha^{-1} u' \mathcal{F}_\alpha(g_\alpha p, s) &=& g_\alpha^{-1} \mathcal{F}_\alpha(g_\alpha p,s)\\
ug_\alpha^{-1} \mathcal{F}_\alpha(g_\alpha p, s) &=& \mathcal{F}_\alpha(p,s)\\
\mathcal{F}_\alpha(up,s) &=& \mathcal{F}_\alpha(p,s)\in g_\alpha^{-1}\A
\end{eqnarray*}
\qed

\begin{lemma}
\label{lem:fillings}
Let $\ell$ be a loop in $U_\alpha(\Z[t])\cdot g_\alpha^{-1}\A$, then there exists an $S>0$ such that $\mathcal{F}_\alpha(\ell,s)\in g_\alpha^{-1}\A$ for all $s\ge S$.
\end{lemma}
\Proof For each point $p\in\ell$, the above lemma gives such an $S$ which we will call $S_p$. As $\ell$ is compact, we have that $\max_{p\in\ell}\{S_p\}$ exists. Let $S$ be this maximum.\qed

\begin{lemma}
Let $\ell$ be a loop in $U_\alpha(\Z[t])\cdot g_\alpha^{-1}\A$, then there exists a filling disk for $\ell$ in $U_\alpha(\Z[t])\cdot g_\alpha^{-1}\A$.
\end{lemma}
\Proof By the above lemma, we have an $S>0$ such that $\mathcal{F}_\alpha(\ell,S)\subseteq g_\alpha^{-1}\A$. Consider the restriction of $\mathcal{F}_\alpha$ to $(U_\alpha(\Z[t])\cdot g_\alpha^{-1}\A)\times[0,S]\to U_\alpha(\Z[t])\cdot g_\alpha^{-1}\A$. This is still continuous and hence is a homotopy. So $\mathcal{F}_\alpha(\ell,[0,S])$ is an annulus in $U_\alpha(\Z[t])\cdot g_\alpha^{-1}\A$ with one boundary component $\ell$ and the other a loop in $g_\alpha^{-1}\A$. Well, $g_\alpha^{-1}\A$ is an apartment, and hence a Euclidean plane. Any loop in a Euclidean plane has a filling disk, so we can attach this filling disk to the annulus, and yield the desired filling disk in $U_\alpha(\Z[t])\cdot g_\alpha^{-1}\A$ for $\ell$.\qed

\subsection{The Global Disk}
In \SecRef{localCocycles}, we detailed a cycle in $\Omega$, the open hemisphere complex of $\Xi$, which we called $\mathfrak{B}$. This cycle was important for the key fact that $\varphi(\pi(\mathfrak{B}))\neq0$. That is, we used this cycle to show that our constructed cocycle was not the zero cycle. In this subsection, we locate this cycle about $y_k$, for each $k>0$, and extend it to be a 2-disk in $X$ which will serve the same purpose for each $\Phi_k$: that it is not the zero cocycle. Furthermore, we do this in such a way that these cycles serve to show the $\Phi_k$ are a linearly independent set.

In \SecRef{localStars}, we defined a map $\eta_y\colon\CH(\St(y))\to\Xi$ for each vertex $y$ in the interior of $\Sect$.  For each $k>0$, let $\eta_k$ be this map for the vertex $y_k$. Define the map $\rho_k$ in a similar fashion. Let $\mathfrak{B}_k = \eta_y^{-1}(\mathfrak{B})$. As $\mathfrak{B}$ is a cycle and $\eta_y$ is an identification, $\mathfrak{B}_k$ must be a disk with boundary a cycle identical to $\mathfrak{B}$ in $\Lk(y_k)$. Furthermore, $\mathfrak{B}_k\subseteq\Omega(y_k)$ from the definition of $\Omega(y_k)$.

\begin{lemma}
For each $k>0$, there exists a disk $\widehat{\mathfrak{B}_k}< X$ such that
\begin{enumerate}
\item $\widehat{\mathfrak{B}_k}\cap (\SL_3(\Q[t])\cdot\St(y_k)) = \mathfrak{B}_k$
\item For each $j>k$, $\widehat{\mathfrak{B}_k}\cap (\SL_3(\Q[t])\cdot\St(y_j)) = \varnothing$
\item $\partial\widehat{\mathfrak{B}_k}\subseteq (U(\Z[t])\cdot\partial\Sect)$.
\end{enumerate}
\end{lemma}

\Proof Let $C$ be a chamber in $\CH(\St(y_k)\cap\A)$. Then we can take the sector in $\A$ with base point $y_k$ such that it includes only $C$ in $\St(y_k)$ and take the intersection with $\Sect$ and call this subset $\widehat{C}$. This is the same construction as in \SecRef{sequenceOfVertices}. Lemma~\ref{lem:integerImage} and the construction of $\mathfrak{B}$ tells us that every chamber $D$ in $\mathfrak{B}_k$ can be written as $D=\gamma C$ where $C\in\CH(\St(y_k)\cap\A)$ and $\gamma\in\Stab_{U(\Z[t])}(y_k)$. For any such chamber $D$, let $\widehat{D}=\gamma\widehat{C}$. Now we can construct
\[
\widehat{\mathfrak{B}_k} = \bigcup_{D \in\CH(\mathfrak{B}_k)} \widehat{D}
\]
It remains to show that this is in fact a compact disk. In \SecRef{sequenceOfVertices}, we showed that each $\widehat{C}$ is compact, and hence each $\widehat{D}$ also is. We are taking a finite union, so $\widehat{\mathfrak{B}_k}$ is also compact. $\widehat{\mathfrak{B}_k}$ is connected as each $\widehat{D}$ includes the vertex $y_k$. We can choose the $\gamma\in\Gamma$ for each chamber in such a way that the only boundary components that occur are contained in $U(\Z[t])(\partial\Sect)$. Notice that in extending each chamber of $\mathfrak{B}_k$, we never introduced any $\SL_3(\Q[t])$-translates of $y_k$, as we stayed within a single translate of $\Sect$. So, when we take the intersection of $\widehat{\mathfrak{B}}_k$ and $\SL_3(\Q[t])\cdot\St(y_k)$, we are left only with $\mathfrak{B}_k$.

The same reasoning gives us part 2 of the lemma. In particular, we specified in our choice of the $y_k$ that each $y_j$ for $j>k$ must not lie in the $\widehat{C}$ used to construct $\widehat{\mathfrak{B}_k}$. Again, as $\Sect$ is a strict fundamental domain for the action of $\SL_3(\Q[t])$ on the building, we are guaranteed to avoid any $\SL_3(\Q[t])$ translates of $y_j$ as well.\qed\\

Now, for each $k>0$, let $L_k=\partial\widehat{\mathfrak{B}_k}\subseteq U(\Z[t])\cdot\partial\Sect$. We can split $L_k$ up into segments which have endpoints in $U(\Z[t])\cdot x_0$ and are each contained completely in $U(\Z[t])\cdot R_\alpha$ for some $\alpha\in\Delta$. That is, there exist paths $L_k^1,\dots, L_k^m$ such that
\begin{enumerate}
\item $L_k^i\cap L^{i+1}_k = u_k^ix_0$ for some $u_k^i\in U(\Z[t])$,
\item For each $i$, there exists an $\alpha\in\Delta$ such that $L^i_k\subseteq U(\Z[t])\cdot R_\alpha$, and
\item $\bigcup_i L_k^i = L_k$
\end{enumerate}

\begin{lemma}
For each $L_k^i$, there exists a disk $D_k^i\subseteq(U_\alpha(\Z[t])\cdot g_\alpha^{-1}\A)$ such that $\partial D_k^i = L_k^i\cup\ell_k^i$, where $\ell_k^i$ is a path in $U_\alpha(\Z[t])\cdot A_\alpha$. Furthermore, $\ell_k^i$ and $L_k^i$ intersect only at their endpoints.
\end{lemma}
\Proof First we remark that $(u_k^i)^{-1} L_k^i$ is a path in $U(\Z[t])\cdot R_\alpha$, and as it contains $x_0$ (and is inherently connected), it is in fact contained in $U_\alpha(\Z[t]) R_\alpha$. Thus, the endpoints of $(u_k^i)^{-1} L_k^i$ are both contained in $U_\alpha(\Z[t])\cdot A_\alpha$. Lemma~\ref{lem:UAalphaConnected} shows that this is a connected space, and hence there exists a path $(u_k^i)^{-1}\ell_k^i\subseteq U_\alpha(\Z[t])\cdot A_\alpha$ which connects these endpoints. So by lemma ~\ref{lem:fillings}, there exists a filling disk in $U_\alpha(\Z[t])\cdot g_\alpha^{-1}\A$ for $(u_k^i)^{-1}\ell_k^i\cup (u_k^i)^{-1}L_k^i$. Let $D_k^i$ be the $u_k^i$ translate of this disk, it satisfies all the criteria for the lemma.\qed

Now we put these pieces together and construct our global disks.

\begin{proposition}
\label{prop:globalDisk}
For each $k>0$, there exists a disk $\mathfrak{D}_k \subseteq X$ such that
\begin{enumerate}
\item $\mathfrak{D}_k\cap (\SL_3(\Q[t])\cdot\St(y_k)) = \mathfrak{B}_k$
\item For $j>k$,  $\mathfrak{D}_k\cap (\SL_3(\Q[t])\cdot\St(y_j)) = \varnothing$
\item $\partial\mathfrak{D}_k\subseteq X_0^{(1)}$. Where $X_0^{(1)}$ is the 1-skeleton of $X_0=\SL_3(\Z[t])\cdot\C$.
\end{enumerate}
\end{proposition}

\Proof Let $\widehat{\mathfrak{D}}_k = \widehat{\mathfrak{B}}_k\cup\bigcup_i D_k^i$. By construction, this is a disk with boundary equal to $\bigcup_i\ell_k^i\subseteq U(\Z[t])\cdot (A_{\alpha_{12}}\cup A_{\alpha_{23}})\subseteq X_0$. So by the cellular mapping theorem, we can form a homotopy of $\partial \widehat{\mathfrak{D}}_k$ to a loop in $X_0$ which is contained in the 1-skeleton. After glueing the annulus given by this homotopy onto $\widehat{\mathfrak{D}}_k$, we end up with a disk $\mathfrak{D}_k$ with boundary in the 1-skeleton of $X_0$.

The disk $\widehat{\mathfrak{B}_k}$ satisfies property 1 and 2. Each disk $D_k^i$ lies within $\SL_3(\Q[t])\cdot g_\alpha\A$ for some $\alpha\in\Delta$. Each $g_\alpha\A$ has fundamental domain under the action of $\SL_3(\Q[t])$ which is contained in the closed star about the boundary of $\Sect$. We chose our $y_k$ specifically so that they do not lie in this region. Hence, none of the $D_k^i$ which are glued in to construct $\mathfrak{D}_k$ will pose a problem for criteria 1 and 2 in the proposition.\qed\\

\subsection{Turning Global Disks Into Quotient Cycles}
In this section, we take our global disks $\mathfrak{D}_k$ that we just constructed and use them to build 2-cycles in $\Gamma\setminus Y$.

Each $\partial\mathfrak{D}_k$ is a loop in $X_0$, and hence we glued in a disk $\mathfrak{d}_k$ during our construction of $Y$ with boundary $\partial\mathfrak{D}_k$. As $X$ is contractible, filling disks are unique. Hence we have that $\psi(\mathfrak{d}_k) = \mathfrak{D}_k$. $\Gamma\mathfrak{d}_k$ is a 2-cell in $\Gamma\setminus Y$ with boundary in $\Gamma X_0$. As this is a contractible space, we know there exists some filling disk $\mathfrak{z}_k$ for $\partial\Gamma\mathfrak{d}_k$ in $\Gamma X_0$. Let:
\[
\mathfrak{h}_k = \mathfrak{z}_k + \Gamma\mathfrak{d}_k
\]

{\it Proof of Proposition~\ref{prop:cocyclesExist}:} We defined the $\mathfrak{h}_k$ just above. First we show that each is evaluated to be non-zero by $\Phi_k$.\\
\begin{eqnarray*}
\Phi_k(\mathfrak{h}_k) &=& \Phi_k(\Gamma\mathfrak{d}_k + \mathfrak{z}_k)\\
&=& \Phi_k(\Gamma\mathfrak{d}_k) + 0\\
&=& \sum_{\gamma U(\Gamma)\in \Gamma / U(\Gamma)} \varphi_k(\gamma^{-1} \mathfrak{d}_k)\\
&=& \sum_{\gamma U(\Gamma)\in \Gamma / U(\Gamma)}\varphi(H_k(\psi(\mathfrak{d}_k)))\\
&=& \sum_{\gamma U(\Gamma)\in \Gamma / U(\Gamma)}\varphi(H_k(\mathfrak{D}_k))
\end{eqnarray*}
where $\Phi_k(\Gamma\mathfrak{z}_k) = 0$ as it is constructed of cells in $X_0$, which are taken to zero.

Now let $\gamma\in\Gamma / U(\Gamma)$. Then $\gamma y_k\neq y_k$, as the $\Gamma$ stabilizer of $y_k$ is contained in $U(\Gamma)$. So by proposition~\ref{prop:globalDisk}, we have that $\gamma^{-1}\mathfrak{D}_k\cap\St(y_k) = \varnothing$, which means that $H_k(\gamma^{-1}\mathfrak{D}_k)=\varnothing$. We proceed,
\begin{eqnarray*}
\Phi_k(\mathfrak{h}_k)  &=& \sum_{\gamma U(\Gamma)\in \Gamma / U(\Gamma)}\varphi(\bar{\eta}_k(\mathfrak{D}_k))\\
&=& \varphi(H_k(\mathfrak{D}_k))\\
&=& \varphi(\pi_k(\mathfrak{B}_k)\\
&=& \varphi(\pi(\mathfrak{B})) = \pm1 \neq 0
\end{eqnarray*}
Now let $j>k$. Proposition~\ref{prop:globalDisk} tells us that no $\Gamma$-translate of $\mathfrak{D}_k$ will every intersect with a $\SL_3(\Q[t])$-translate of $y_j$. This means that $\Phi_k(\mathfrak{h}_k)=0$. 
\qed\\

{\it Proof of Proposition~\ref{prop:infiniteDimensional}:} We have just constructed an infinite family of 2-cycles $\{\mathfrak{h}_k\}$ in $H_2(\Gamma\setminus Y; \Q)$ . Proposition~\ref{prop:cocyclesExist} shows us that the each $\mathfrak{h}_j$ is not in the span of $\{\mathfrak{h}_k\}_{k=1}^{j-1}$, so $H_2(\Gamma\setminus Y; \Q)$ must be infinite dimensional.\qed

\section{Proof of the Main Theorem}
\label{sec:mainProof}

\begin{lemma}
$H_i(\Gamma, C(Y;\Q))\cong H_i(\Gamma; \Q)$ for $i < 3$.
\end{lemma}
\Proof This is a consequence of $Y$ being $2$-connected and proposition 7.3 in \cite{Br}.\qed\\

Now we have a spectral sequence as described in (7.7) of \cite{Br} which has:
\[
E_{pq}^1 = \bigoplus_{\sigma\in \Sigma_p} H_q(\Gamma_\sigma; \Q) \implies  H_{p+q}(\Gamma, C(Y, \Q)) = H_{p+q}(\Gamma; \Q)
\]
where $\Sigma_p$ indexes over the $p$-cells in $\Gamma\setminus Y$. That is:
\[
H_{2}(\Gamma; \Q) = \bigoplus_{p+q=2}E^\infty_{p,q}
\]
We will show that $E_{2,0}^\infty$ is infinite dimensional, which will prove the theorem. Recall that $E_{2,0}^\infty = E_{2,0}^r$ for sufficiently large $r$, so we will show that $E_{2,0}^r$ is infinite dimensional for all $r\ge2$.

\begin{proposition}
\label{prop:infE}
$E_{2,0}^r$ is infinite dimensional for $r\ge 2$.
\end{proposition}

First we prove the following three lemmas:
\begin{lemma}
\label{lem:finiteHomology}
Let $\sigma < X$ be a cell, then each homology group of $\Gamma_\sigma$ is finite dimensional.
\end{lemma}
\Proof Bux-Mohammadi-Wortman proved that $\SL_n(\Z[t])_\sigma$ is $FP_\infty$ for any such $\sigma$ as lemma 2 of \cite{BMW}. As $\Gamma$ is finite index in $\G$, $\Gamma_\sigma$ is also finite index in $\SL_n(\Z[t])_\sigma$. It is well known that being $FP_\infty$ is fixed under finite index subgroups (see proposition 5.1 in \cite{Br}), thus $\Gamma_\sigma$ is also $FP_\infty$. This implies that $\Gamma_\sigma$ has finitely generated homology in every dimension.\qed

\begin{lemma}
\label{lem:finDim}
If $r,q\ge1$, then $E_{p,q}^r$ is finite dimensional.
\end{lemma}
\Proof As $[\SL_3(\Z[t])\ \colon\ \Gamma]<\infty$, $\Gamma$ acts cocompactly on $X_0 = \SL_3(\Z[t])\cdot\C$, and freely on $Y\setminus X_{0}$ by construction. Thus there are only finitely many $c\in \Sigma_p$ such that $\Gamma_c\neq 1$. So $E_{pq}^1 = \oplus_{c\in\Sigma_p} H_q(\Gamma_c; \Q)$ has only finitely many nontrivial summands, and by lemma~\ref{lem:finiteHomology}, each is finitely generated. Thus $E_{pq}^1$ itself is finitely generated.

Recall by definition that $E_{pq}^{r+1} = \ker(d^r_{pq}) / \im(d^r_{p+r,q-r+1})$, where:
\begin{eqnarray*}
d^r_{pq}\colon E_{pq}^r &\to& E_{p-r, q+r-1}^r\\
d^r_{p+r,q-r+1}\colon E_{p+r,q-r+1}^r &\to& E_{pq}^r\\
\end{eqnarray*}
Thus the dimension of $E^{r+1}_{pq}$ is bounded by the dimension of $E^r_{pq}$, proving the lemma.\qed

The proof of the lemma below follows lemma 21 in \cite{Wort}.
\begin{lemma}
\label{lem:equivalence}
$E_{p,0}^2 = H_p(\Gamma\setminus Y; \Q)$.
\end{lemma}
\Proof Let $\partial'$ be the boundary operator for $C_*(Y; \Q)$ and for any $(p-1)$-cell $d$ in $Y$, let:
\begin{eqnarray*}
\pi_d\colon C_{p-1}(Y;\Q) &\to& C_{p-1}(Y; \Q)\\
\sum_{c < Y; \dim(c) = p-1} a_c\cdot c &\mapsto& a_d\cdot d
\end{eqnarray*}
that is, $\pi_d$ is projection.

Let $\partial$ be the boundary operator for the chain complex $C_*(\Gamma\setminus Y; \Q)$. Using the same definition as in the proof of the lemma 10.4, we see that $E_{*,0}^2$ is the homology of the complex $(E_{k,0}^1, d^1)$. We can write:
\[
E_{p,0}^1 = \bigoplus_{c\in\Sigma_p} H_0(\Gamma_c; \Q) = \bigoplus_{c\in\Sigma_p} \Q
\]
and identify this module with $C_p(\Gamma\setminus Y; \Q)$ by the map:
\[
(a_c)_{c\in\Sigma_p} \mapsto \sum_{\Gamma c\subseteq\Gamma\setminus Y} a_c(\Gamma c)
\]
where each $a_c\in\Q$. So in order to prove the claim, it suffices to show that the map:
\[
d^1_{p,0} \colon E_{p,0}^1 = C_p(\Gamma\setminus Y; \Q) \to E_{p-1,0}^1 = C_{p-1}(\Gamma\setminus Y; \Q)
\]
is equivalent to $\partial$ for each $p$.

Let $D_c$ be the set of $(p-1)$-cells in $Y$ contained in $c$. Following section 8 of \cite{Br}, we define the following three homomorphisms:
\[
t_c\colon H_0(\Gamma_c; \Q) \to H_0(\Gamma_{c'}\cap\Gamma_c; \Q)=H_0(\Gamma_c; \Q)
\]
for any $c\in\Sigma_p$, and $c'\in D_c$ to be the transfer map. But as $H_0(\Gamma_c; \Q) = \Q$, $t^c=id$ for any cell $c$. And:\\
\[
u_{cc'}\colon H_0(\Gamma_c; \Q) \to H_0(\Gamma_{c'}; \Q)
\]
which is induced by $\Gamma_c\hookrightarrow\Gamma_{c'}$ and $\pi_d\circ\partial'|_c$. It is clear that:
\[
\sum_{c'\in D_c} u_{cc'}(a_c) = \partial(a_c(\Gamma c))
\]
where we have identified $a_c\in H_0(\Gamma_c; \Q) = \Q$ with $a_c\in\Q$.\\
Lastly:
\[
v_{c'}\colon H_0(\Gamma_{c'}; \Q) \to H_0(\Gamma_{c'}; \Q)
\]
which is induced by conjugation. As $\Q$ is abelian, this is also the identity.

\cite{Br} tells us that, up to sign, the map:
\begin{eqnarray*}
\lambda\colon \bigoplus_{c\in\Sigma_p} H_0(\Gamma_c; \Q) &\to& \bigoplus_{c'\in\Sigma_{p-1}} H_0(\Gamma_{c'}; \Q)\\
\lambda | H_0(\Gamma_c; \Q) &=& \sum_{c'\in D_c} v_{c'}u_{cc'}t_c
\end{eqnarray*}
is the map $d^1$. So by above, we have in fact shown that up to sign, $\partial = d^1$. This proves the claim.\qed\\

{\it{Proof of Proposition~\ref{prop:infE}.}} We proceed by induction on $r$.\\
For the base case, when $r=2$, lemma~\ref{lem:equivalence} tells us that $E_{2,0}^2 = H_2(\Gamma\setminus Y; \Q)$, which we showed is infinite dimensional as proposition~\ref{prop:infiniteDimensional}.

For the inductive case, we assume that $E_{2,0}^r$ is infinite dimensional, and $E_{2,0}^{r+1}$ is the homology of the sequence:
\[
E_{2+r, 1-r}^r \xrightarrow{d^r_{2+r, 1-r}} E_{2,0}^r \xrightarrow{d^r_{2,0}} E_{2-r, r-1}^r
\]
As the spectral sequence is 0 outside the first quadrant, we know that in fact, $E_{2+r, 1-r}^r=0$, and lemma~\ref{lem:finDim} tells us that $E_{2-r, r-1}^r$ is finite dimensional. Thus the kernel of $d^r_{2,0}$ must be infinite dimensional and the image of $d^r_{2+r, 1-r}$ is 0. Thus $E_{2,0}^{r+1}$ is also infinite dimensional.\qed\\

{\it{Proof of Theorem 1:} }Recall that:
\[
H_{2}(\Gamma; \Q) = \bigoplus_{p+q=2}E^\infty_{p,q}
\]
Proposition~\ref{prop:infE} tells us that $E^\infty_{2,0}$ is infinite dimensional. Thus $\bigoplus_{p+q=2}E^\infty_{p,q}=H_{2}(\Gamma; \Q)$ must be infinite dimensional. This tells us that $H^{2}(\Gamma; \Q)$ is infinite dimensional as well.\qed

\bibliography{GoroffThesis}

\begin{thebibliography}{10}
\ifx \showADDRESS      \undefined \def \showADDRESS      #1{#1}          \fi
\ifx \showARTICLENO    \undefined \def \showARTICLENO    #1{Article #1}  \fi
\ifx \showAUTHOR       \undefined \def \showAUTHOR       #1{#1}          \fi
\ifx \showAUTHORRAW    \undefined \def \showAUTHORRAW    #1{}            \fi
\ifx \showBIBTYPE      \undefined \def \showBIBTYPE      #1#2{}          \fi
\ifx \showBOOKPAGES    \undefined \def \showBOOKPAGES    #1{#1}          \fi
\ifx \showBOOKTITLE    \undefined \def \showBOOKTITLE    #1{#1}          \fi
\ifx \showCHAPTER      \undefined \def \showCHAPTER      #1{#1}          \fi
\ifx \showCODEN        \undefined \def \showCODEN        #1{CODEN #1}    \fi
\ifx \showCROSSREF     \undefined \def \showCROSSREF     #1{#1}          \fi
\ifx \showDAY          \undefined \def \showDAY          #1{#1}          \fi
\ifx \showDOI          \undefined \def \showDOI        { doi:\penalty 0} \fi
\ifx \url              \undefined \input path.sty \let \url = \path      \fi
\ifx \showEDITION      \undefined \def \showEDITION      #1{#1}          \fi
\ifx \showEDITOR       \undefined \def \showEDITOR       #1{#1}          \fi
\ifx \showEDITORRAW    \undefined \def \showEDITORRAW    #1{}            \fi
\ifx \showHOWPUBLISHED \undefined \def \showHOWPUBLISHED #1{#1}          \fi
\ifx \showINSTITUTION  \undefined \def \showINSTITUTION  #1{#1}          \fi
\ifx \showISBN         \undefined \def \showISBN         #1{ISBN #1}     \fi
\ifx \showISBNXIII     \undefined \def \showISBNXIII     #1{ISBN-13 #1}  \fi
\ifx \showISSN         \undefined \def \showISSN         #1{ISSN #1}     \fi
\ifx \showISSNL        \undefined \def \showISSNL        #1{ISSN-L #1}   \fi
\ifx \showJOURNAL      \undefined \def \showJOURNAL      #1{#1}          \fi
\ifx \showKEY          \undefined \def \showKEY          #1{#1}          \fi
\ifx \showLCCN         \undefined \def \showLCCN         #1{LCCN #1}     \fi
\ifx \showMONTH        \undefined \def \showMONTH        #1{#1}          \fi
\ifx \showNOTE         \undefined \def \showNOTE         #1{#1}          \fi
\ifx \showNUMBER       \undefined \def \showNUMBER       #1{#1}          \fi
\ifx \showORGANIZATION \undefined \def \showORGANIZATION #1{#1}          \fi
\ifx \showPAGECOUNTONE \undefined \def \showPAGECOUNTONE #1{, #1 page}   \fi
\ifx \showPAGECOUNT    \undefined \def \showPAGECOUNT    #1{, #1 pages}  \fi
\ifx \showPAGES        \undefined \def \showPAGES        #1{#1}          \fi
\ifx \showPRICE        \undefined \def \showPRICE        #1{#1}          \fi
\ifx \showPUBLISHER    \undefined \def \showPUBLISHER    #1{#1}          \fi
\ifx \showSCHOOL       \undefined \def \showSCHOOL       #1{#1}          \fi
\ifx \showSERIES       \undefined \def \showSERIES       #1{#1}          \fi
\ifx \showTITLE        \undefined \def \showTITLE        #1{#1}          \fi
\ifx \showTYPE         \undefined \def \showTYPE         #1{#1}          \fi
\ifx \showURL          \undefined \def \showURL         {URL }           \fi
\ifx \url              \undefined \input path.sty \let \url = \path      \fi
\ifx \showVOLUME       \undefined \def \showVOLUME       #1{#1}          \fi
\ifx \showYEAR         \undefined \def \showYEAR         #1{#1}          \fi

\bibitem{BoSe}
  \ifshowBIBTYPE \showBIBTYPE{article}{BoSe} \fi \showAUTHORRAW{Borel, A.,
  Serre, J-P.}\showAUTHOR{Serre J-P. Borel, A.}
\newblock \showTITLE{Corners and arithmetic groups}.
\newblock {\em \showJOURNAL{Commentarii mathematici Helvetici}},
  \showVOLUME{48}:\penalty 0 \showPAGES{436--483}, \showYEAR{1973}.
\newblock \ifshowURL {\showURL \url{http://eudml.org/doc/139559}}. \fi

\bibitem{Build}
  \ifshowBIBTYPE \showBIBTYPE{book}{Build} \fi \showAUTHORRAW{Kenneth
  Brown}\showAUTHOR{Kenneth Brown}.
\newblock {\em \showTITLE{Buildings}}.
\newblock \showPUBLISHER{Springer-Verlag New York}, \showYEAR{1989}.

\bibitem{Br}
  \ifshowBIBTYPE \showBIBTYPE{book}{Br} \fi \showAUTHORRAW{Kenneth S.
  Brown}\showAUTHOR{Kenneth~S. Brown}.
\newblock {\em \showTITLE{Cohomology of Groups}}.
\newblock \showSERIES{Graduate Texts in Mathematics}.
  \showPUBLISHER{Springer-Verlag New York}, \showYEAR{1982}.

\bibitem{BKW}
  \ifshowBIBTYPE \showBIBTYPE{article}{BKW} \fi \showAUTHORRAW{Bux, Kai-Uwe and
  K{\"o}hl, Ralf and Witzel, Stefan}\showAUTHOR{Kai-Uwe Bux, Ralf K{\"o}hl, and
  Stefan Witzel}.
\newblock \showTITLE{Higher finiteness properties of reductive arithmetic
  groups in positive characteristic: The rank theorem}.
\newblock {\em \showJOURNAL{Annals of Mathematics}}, \showVOLUME{177}\penalty 0
  (\showNUMBER{1}):\penalty 0 \showPAGES{311--366}, \showMONTH{Jan}
  \showYEAR{2013}. \ifshowISSN {\showISSN{0003-486X}}. \fi
\newblock \ifshowURL {\showURL
  \url{http://dx.doi.org/10.4007/annals.2013.177.1.6}}. \fi
\newblock \ifshowDOI {\showDOI
  \url{10.4007/annals.2013.177.1.6}}\ifshowDOIPERIOD . \fi \fi

\bibitem{BMW}
  \ifshowBIBTYPE \showBIBTYPE{article}{BMW} \fi \showAUTHORRAW{{Bux}, Kai-Uwe
  and {Mohammadi}, Amir and {Wortman}, Kevin}\showAUTHOR{Kai-Uwe {Bux}, Amir
  {Mohammadi}, and Kevin {Wortman}}.
\newblock \showTITLE{{SL(n,Z[t]) is not FP\_\{n-1\}}}.
\newblock {\em \showJOURNAL{arXiv e-prints}}, page \showPAGES{arXiv:0801.1332},
  \showMONTH{Jan} \showYEAR{2008}.

\bibitem{CK}
  \ifshowBIBTYPE \showBIBTYPE{article}{CK} \fi \showAUTHORRAW{Morgan Cesa and
  Brendan Kelly}\showAUTHOR{Morgan Cesa and Brendan Kelly}.
\newblock \showTITLE{{$H^2(SL_3(\mathbb{Z}[t]); \mathbb{Q})$} is infinite
  dimensional}.
\newblock {\em \showJOURNAL{arXiv: Group Theory}}, \showYEAR{2015}.

\bibitem{KM}
  \ifshowBIBTYPE \showBIBTYPE{article}{KM} \fi \showAUTHORRAW{Sava Krsti{\'c}
  and James McCool}\showAUTHOR{Sava Krsti{\'c} and James McCool}.
\newblock \showTITLE{Presenting gl$_n(k\langle t\rangle)$}.
\newblock {\em \showJOURNAL{Journal of Pure and Applied Algebra}},
  \showVOLUME{141}\penalty 0 (\showNUMBER{2}):\penalty 0 \showPAGES{175 --
  183}, \showYEAR{1999}. \ifshowISSN {\showISSN{0022-4049}}. \fi
\newblock \ifshowURL {\showURL
  \url{http://www.sciencedirect.com/science/article/pii/S002240499800022X}}.
  \fi
\newblock \ifshowDOI {\showDOI
  \url{10.1016/S0022-4049(98)00022-X}}\ifshowDOIPERIOD . \fi \fi

\bibitem{Sch}
  \ifshowBIBTYPE \showBIBTYPE{article}{Sch} \fi \showAUTHORRAW{Schulz,
  Bernd}\showAUTHOR{Bernd Schulz}.
\newblock \showTITLE{Spherical subcomplexes of spherical buildings}.
\newblock {\em \showJOURNAL{Geometry and Topolology}}, \showVOLUME{17}\penalty
  0 (\showNUMBER{1}):\penalty 0 \showPAGES{531--562}, \showYEAR{2013}.
\newblock \ifshowURL {\showURL \url{https://doi.org/10.2140/gt.2013.17.531}}.
  \fi
\newblock \ifshowDOI {\showDOI \url{10.2140/gt.2013.17.531}}\ifshowDOIPERIOD .
  \fi \fi

\bibitem{Sou}
  \ifshowBIBTYPE \showBIBTYPE{article}{Sou} \fi \showAUTHORRAW{Christophe
  Soul{\'e}}\showAUTHOR{Christophe Soul{\'e}}.
\newblock \showTITLE{Chevalley groups over polynomial rings}.
\newblock {\em \showJOURNAL{Homological Group Theory}}, pages
  \showPAGES{359--367}, \showYEAR{1977}.

\bibitem{Sus}
  \ifshowBIBTYPE \showBIBTYPE{article}{Sus} \fi \showAUTHORRAW{Suslin,
  A.A.}\showAUTHOR{A.A. Suslin}.
\newblock \showTITLE{On the structure of the special linear group over
  polynomial rings}.
\newblock {\em \showJOURNAL{Mathematics of the {USSR}-Izvestiya}},
  \showVOLUME{11}\penalty 0 (\showNUMBER{2}):\penalty 0 \showPAGES{221--238},
  \showYEAR{1977}.
\newblock \ifshowDOI {\showDOI
  \url{10.1070/IM1977v011n02ABEH001709}}\ifshowDOIPERIOD . \fi \fi

\bibitem{Wort}
  \ifshowBIBTYPE \showBIBTYPE{article}{Wort} \fi \showAUTHORRAW{Wortman,
  Kevin}\showAUTHOR{Kevin Wortman}.
\newblock \showTITLE{An infinitely generated virtual cohomology group for
  noncocompact arithmetic groups over function fields}.
\newblock {\em \showJOURNAL{Groups, Geometry, and Dynamics}},
  \showVOLUME{10}\penalty 0 (\showNUMBER{1}):\penalty 0 \showPAGES{91--115},
  \showYEAR{2016}. \ifshowISSN {\showISSN{1661-7207}}. \fi
\newblock \ifshowURL {\showURL \url{http://dx.doi.org/10.4171/GGD/344}}. \fi
\newblock \ifshowDOI {\showDOI \url{10.4171/ggd/344}}\ifshowDOIPERIOD . \fi \fi

\end{thebibliography}


\providecommand{\bysame}{\leavevmode\hbox to3em{\hrulefill}\thinspace}
\providecommand{\MR}{\relax\ifhmode\unskip\space\fi MR }
\providecommand{\MRhref}[2]{%
  \href{http://www.ams.org/mathscinet-getitem?mr=#1}{#2}
}
\providecommand{\href}[2]{#2}
\bibliographystyle{ieeetr}

\end{document}